\documentclass[12pt]{amsart}
\usepackage{amsmath,amssymb,amsfonts,graphicx,graphics, setspace}

\newcommand{\mathsym}[1]{{}}
\newcommand{\unicode}[1]{{}}

\usepackage{tikz}
\usetikzlibrary{shapes,snakes}

\newtheorem{thm}{Theorem}

\theoremstyle{definition}

\newtheorem{lemma}{Lemma}[section]
\newtheorem{proposition}[lemma]{Proposition}
\newtheorem{condition}{Condition}[section]
\theoremstyle{remark}
\newtheorem{remark}[lemma]{Remark}
\numberwithin{equation}{section}
\newtheorem{cor}{Corollary}[section]

\theoremstyle{definition}
\newtheorem{defi}[cor]{Definition}

\newcommand{\co}{\colon\thinspace}

\begin{document}

\title{``Slicing'' the Hopf link}
\author{Vyacheslav Krushkal}
\address{Department of Mathematics, University of Virginia, Charlottesville, VA 22904}
\email{krushkal\char 64 virginia.edu}

\thanks{This research was supported in part by the NSF}

\begin{abstract} 
A link in the $3$-sphere is called (smoothly) slice if its components bound disjoint smoothly embedded disks in 
the $4$-ball. More generally, given a $4$-manifold $M$ with a distinguished circle in its boundary, a link in the 
$3$-sphere is called {\em $M$-slice} if its components bound in the $4$-ball disjoint embedded copies of $M$.
%The Hopf link is in a sense ``the most non-slice'' link, and the non-trivial linking number provides a
%well-known and elementary obstruction to sliceness. 
A $4$-manifold $M$ is constructed such that  the Borromean rings are not $M$-slice but the Hopf link is. This contrasts
the classical link-slice setting where the Hopf link may be thought of ``the most non-slice'' link.
Further examples and an obstruction for a family of decompositions of the $4$-ball are discussed in the context of the A-B slice problem.
\end{abstract}

\maketitle

\section{Introduction} \label{introduction}
The classification of knots and links up to concordance, and in particular the study of slice links,
is a classical and challenging problem at the interface between $3$- and $4$-manifold 
topology. Recall that a link in the $3$-sphere is called smoothly (respectively topologically) {\em slice} if its components bound disjoint smooth 
(respectively locally flat) embedded disks in the $4$-ball, where $S^3=\partial D^4$. The results of this paper
take place in the smooth category, so without further mention all $4$-manifolds and maps between them will be smooth.

Let $M$ be a $4$-manifold with a distinguished circle in its boundary ${\gamma}\subset \partial M$, embedded in the $4$-ball: $(M, {\gamma})\subset (D^4, S^3)$. A link $L=(l_1\ldots, l_n)$ in the 
$3$-sphere is called {\em $M$-slice} if there exist $n$ disjoint embeddings
$f_i\co (M, {\gamma})\hookrightarrow (D^4, S^3)$ 
such that $f_i({\gamma})=l_i, i=1,\ldots, n$. 
The classical notion of a slice link
corresponds to $M$ equal to the $2$-handle: $(M,{\gamma})=(D^2\times D^2, \partial D^2\times 0)$. The more general notion of $M$-slice links is considerably more subtle, with the topology of $M$ playing an important role. In particular, a new feature not present in the classical setting is that $1$- and $2$-handles of $M$ may link when $M$ is embedded in $D^4$. We impose an additional requirement in the definition of $M$-slice that each embedding $f_i$ is isotopic to the original embedding  $(M, {\gamma})\subset (D^4, S^3)$. This condition, motivated by the A-B slice problem (see section \ref{ABdecomposition subsection}), gives some control over the complexity of the problem, for example allowing one to keep track of the linking of the handles of $M$ in $4$-space. (It follows from
\cite{K4} that the resulting theory is quite different depending on whether this requirement is imposed or not.
This is discussed in more detail further below.) The main result of this paper is the following theorem.

\begin{thm} \label{main theorem} \sl There exist $4$-manifolds $M$ such that the Hopf link is $M$-slice but the Borromean rings are not.
\end{thm}

This result contrasts the usual slice setting: note that any link in $S^3$ bounds smooth disks in $D^4$, possibly intersecting and self-intersecting in a finite
number of transverse double points.
The intersections (and self-intersections) of surfaces in $4$-space are locally modeled on two coordinate planes intersecting at the origin 
in ${\mathbb R}^4$, and the link of the singularity is the Hopf link. Thus in a naive sense, ``if the Hopf link were slice the double points could be resolved
and any link would be slice''.  In this imprecise sense, the Hopf link plays a role in link theory analogous to that of the free group $F$ on two 
generators in group theory: the Hopf link is ``the most non-slice link'' similarly to $F$ being ``the most non-abelian group''.
The theorem above  shows that this analogy does not extend to the more general notion of $M$-slice links.

The manifold $M$ in the proof of the theorem is constructed in section \ref{construction section} as a handlebody with $1$- and $2$-handles, 
and the embedding/non-embedding results are proved in the {\em relative-slice} context introduced in \cite{FL}. Here the
$1$-handles embedded in $D^4$ are dually considered as $2$-handles removed from the collar on 
the boundary, and the embedding question is equivalent to slicing the attaching link for the $2$-handles ``relative to'' the link corresponding to the $1$-handles.
The fact that the Hopf link is $M$-slice is proved in section \ref{Hopf section}, 
using properties of the Milnor group \cite{M} adapted to the relative-slice context. 
The second part of the theorem, asserting that the Borromean rings are not
$M$-slice, relies on a subtle calculation in commutator calculus; see section
\ref{Borromean section}.

The notion of an M-slice link arises naturally in the the {\em A-B slice problem} \cite{F2, FL}, a formulation of 
the $4$-dimensional topological surgery conjecture.
The analysis necessary
for finding an obstruction for the Borromean rings in theorem \ref{main theorem} is quite different
compared to previously considered examples in the subject.
In particular, this is the first observed case
where the system of equations associated to the relative-slice problem has
{\em rational} but not integral solutions.
If the answer to the A-B slice problem turns out to be positive (i.e. if surgery works for free groups), then it seems likely that the phenomena observed in this
paper will have a role in constructing the relevant A-B decompositions.
On the other hand, the results of the paper are consistent with the conjecture \cite{F2}
that the Borromean rings are not A-B slice,  see section \ref{AB slice section} for further discussion.

\section{Construction of $M$} \label{construction section}

\begin{figure}[h]
\includegraphics[height=3.7cm]{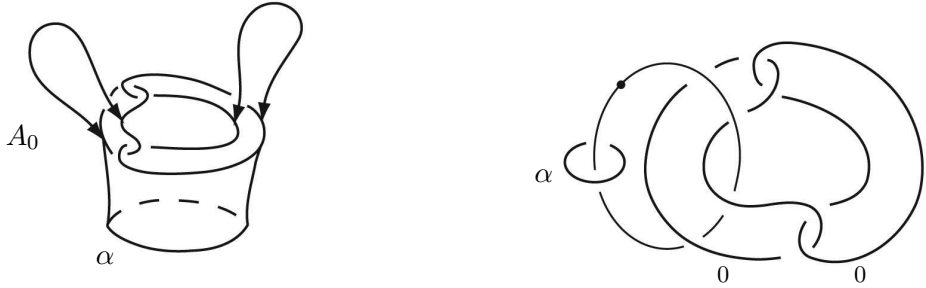} 
{\small
    \put(-316,5){$\alpha$}
    \put(-350,50){$A_0$}
    \put(-150,36){$\alpha$}}
{\scriptsize
    \put(-81,-2){$0$}
    \put(-29,-2){$0$}}
\caption{A preliminary construction: zero-framed $2$-handles
attached to the Bing double, a schematic picture of the spine and a Kirby diagram.}
\label{fig:Bing pair}
\end{figure}

The starting point of the construction is the handlebody $A_0=S^1\times D^2\times [0,1]\, \cup\, $(two zero-framed 
$2$-handles), where the $2$-handles are attached to the Bing double of the core of the solid torus
$S^1\times D^2\times \{1\}$, figure \ref{fig:Bing pair}.  

The distinguished circle (``attaching curve'') of $A_0$
is ${\alpha}=S^1\times \{0\}\times\{0\}$. 
This handlebody is easily seen to embed into $D^4$; it is 
the complement of a standard embedding into $D^4$ of a genus one surface with boundary: $D^4=A_0\cup B_0$
where $B_0={\Sigma}\times D^2$, $\Sigma$ is a genus surface with $\partial {\Sigma}={\beta}$, and the curves
${\alpha}, {\beta}$ form the Hopf link in $\partial D^4$.
Iterating this construction (applying the Bing doubling described above) to various $2$-handles
of $A_0, B_0$ one gets the family of {\em model decompositions} of the $4$-ball, see \cite{FL} and also \cite{FK}
for more details. The construction in this paper builds on recent work of the author in \cite{K4, K5}, and it
is quite different from the model decompositions. 

\begin{figure}[ht]
%\centering
%\vspace{.2cm}
\includegraphics[height=3.7cm]{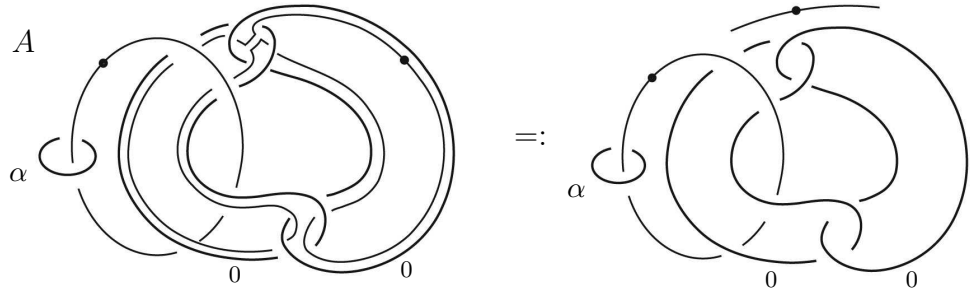}
{\small
    \put(-366,36){$\alpha$}
    \put(-155,30){$\alpha$}}
    \put(-365,85){$A$}
    \put(-175,50){$=:$}
{\scriptsize
    \put(-283,-2){$0$}
    \put(-218,0){$0$}
    \put(-80,-4){$0$}
    \put(-27,-4){$0$}}
%    \vspace{.45cm} 
\caption{A Kirby diagram of the $4$-manifold $A$. 
The figure on the right will serve as a short-hand notation for $A$ in follow-up sections to avoid drawing the complicated dotted curve.}
\label{fig:Bnew}
\end{figure}

The relevant $4$-manifold $A$ used in the
proof of theorem \ref{main theorem} is obtained by attaching a single $1$-handle to $A_0$, as shown on the left in figure
\ref{fig:Bnew}. (The actual manifold $M$ in the statement of theorem \ref{main theorem} will be defined as $A$ with a number of self-plumbings 
of its $2$-handles, see section \ref{completion}.) To avoid drawing unnecessarily complicated diagrams later in the paper, 
a short hand handle notation for $A$ is introduced on the right in figure \ref{fig:Bnew}.

Note that the link formed by the two dotted components is the two-component unlink, and the diagram in figure \ref{fig:Bnew} is indeed
a Kirby diagram of a $4$-manifold. There is a band visible in the picture which is involved in the connected sum of ``parallel copies'' of the two 
zero-framed $2$-handles (they are not actual parallel copies since the dotted curve on the left and the attaching curves of the two $2$-handles form the
Borromean rings, while the dotted curve on the left and the two ``parallel copies'' form the unlink as shown on the right in figure \ref{fig:more examples}). A particular choice of this band in
the $3$-sphere is not going to be important for the argument, as long as the dotted link is the unlink. The construction in figure
\ref{fig:Bnew} differs from an example in \cite{K4, K5} in the attaching curve of the ``interesting'' $1$-handle. The properties
of the resulting $4$-manifolds are quite different, and the analysis required to formulate an obstruction for the Borromean rings
in this paper is substantially more subtle. This work sheds a new light on the techniques necessary for a solution to the A-B slice problem,
see section \ref{AB slice section}.

\section{The Hopf link is $M$-slice} \label{Hopf section}

The proof of theorem \ref{main theorem} will be given in the following two sections in
the context of the {\em relative-slice} problem
(introduced in \cite{FL} and also described below) using the Milnor group. The reader is referred
to the original reference \cite{M} for a more complete introduction to the Milnor
group of links in the $3$-sphere. The application in this paper will concern a variation of the theory for 
submanifolds in $4$-space which will be summarized next.

\subsection{The Milnor group.} \label{Milnor group subsection}

\begin{defi}
Let $G$ be a group normally generated by a finite collection of elements $g_1,\ldots, g_n$. The 
{\em Milnor group} of $G$, relative to the given normal generating set $\{ g_i\}$,
is defined as
\begin{equation} \label{eq:Milnor group}
MG := G \, /\, \langle\! \langle \, [g_i^x, g_i^y] \;\; i=1,\ldots, n,\;\, x,y\in G \rangle\!\rangle.
\end{equation}
\end{defi}

The Milnor group is a finitely presented nilpotent group of class $\leq n$, where
$n$ is the number of normal generators in the definition above, see \cite{M}.
Suppose ${\Sigma}$ is a collection of surfaces with boundary, properly and disjointly embedded in
$(D^4, S^3)$, and let $G$ denote ${\pi}_1(D^4\smallsetminus {\Sigma})$.

Consider meridians
$m_i$ to the components ${\Sigma}_i$ of ${\Sigma}$: $m_i$ is an
element of $G$ which is obtained by following a path $\alpha_i$ in $D^4\smallsetminus {\Sigma}$
from the basepoint to the boundary of a regular neighborhood of ${\Sigma}_i$, followed
by a small circle (a fiber of the circle normal bundle) linking ${\Sigma}_i$, then followed
by ${\alpha}_i^{-1}$. Observe that $G$ is normally generated by the elements $\{ m_i\}$, one for each
component of $\Sigma$.

Let $F_{g_1,\ldots, g_n}$ denote the free group generated by the $\{g_i\}$, $i=1,\ldots, n$, and consider the Magnus
expansion
\begin{equation}\label{Magnus}
M\co F_{g_1,\ldots, g_n}\longrightarrow {\mathbb Z}[\! [x_1,\ldots, x_n]\! ]
\end{equation}
into the ring
of formal power series in non-commuting variables $\{ x_i\}$, defined by $$M(g_i)=1+x_i, \;\,
M(g_i^{-1})=1-x_i+x_i^2-x_i^3\pm\ldots$$ The Magnus expansion induces a homomorphism (which abusing the notation we denote
again by $M$) from the free Milnor group
\begin{equation} \label{MagnusMilnor}
M\co MF_{g_1,\ldots, g_n}\longrightarrow R_{x_1,\ldots,x_n},
\end{equation}
into the quotient $R_{x_1,\ldots,x_n}$ of ${\mathbb Z}[\! [x_1,\ldots, x_n]\!  ]$ by the ideal generated by all monomials
$x_{i_1}\cdots x_{i_k}$ with some index occuring at least twice. It is proved in \cite{M} that
the homomorphism (\ref{MagnusMilnor}) is well-defined and injective.

The relations in (\ref{eq:Milnor group}) are very well suited for studying 
links $L$ in $S^3$ up to {\em link homotopy}. In this original setting for the definition
of the Milnor group \cite{M} one takes $G$ to be the link group
${\pi}_1(S^3\smallsetminus L)$ and a normal set of generators is provided by
meridians $m_i$  to the link components. Two links are link-homotopic if they
are connected by a $1$-parameter family of link maps where different components
stay disjoint for all values of the parameter. If $L$, $L'$ are link-homotopic then 
their Milnor groups $ML$, $ML'$ are isomorphic, and moreover an $n$-component link $L$ is
null-homotopic in this sense if and only if $ML$ is isomorphic to the free Milnor group $MF_{m_1,\ldots, m_n}$.

The Milnor group is also useful for studying surfaces $\Sigma$ in the $4$-ball which are disjoint
but which may have self-intersections: in this case the Clifford tori linking
the double points in $D^4$ give rise to the relations (\ref{eq:Milnor group}) in $M{\pi}_1(D^4\smallsetminus {\Sigma})$.
The theory of link homotopy discussed above may be interpreted as the study of links up to singular concordance
(links $L\subset S^3\times\{0\}$, $L'\subset S^3\times\{1\}$ bounding disjoint maps of annuli into $S^3\times[0,1]$).
In particular, a link is null-homotopic if and only if its components bound disjoint maps of disks into $D^4$, and in this case the Milnor
group is isomorphic to the free Milnor group. Surfaces of higher genus give rise to  
additional relations in the fundamental group of the complement, and the Milnor group in this more general case depends on
linking of the surfaces in $D^4$.

\subsection{The relative slice problem.} \label{relative slice subsection}

The proofs of the embedding and non-embedding statements in this paper will be based on the
{\em relative-slice} reformulation of the problem, see \cite{FL} and also \cite{K4} for
a more detailed introduction. 

Let $L=\{ l_i\}$ denote the attaching curves of the $2$-handles
of the $4$-manifolds that are to be embedded in the $4$-ball, also let $R=\{ r_j\}$ denote the
dotted curves corresponding to the $1$-handles. 
The $1$-handles are considered as unknotted $2$-handles removed
from the collar on the attaching region of a given $4$-manifold. Considering a slightly smaller $4$-ball
(the original $D^4$ minus the collars on the attaching regions), a given embedding problem is then equivalent to ``slicing $L$ relative to $R$'':
finding slices for the link $L$ in the handlebody $D^4\cup_R 2$-handles, where the $2$-handles are attached with
zero framing to $D^4$ along the dotted components $R$. (See \cite{FL, K4} for more details and illustrations.)

\begin{figure}[ht]
\includegraphics[height=3.7cm]{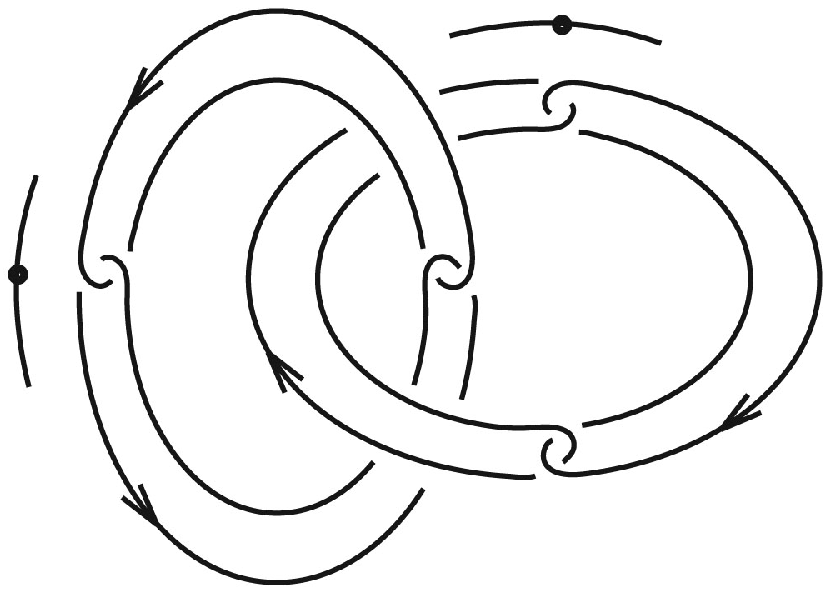}
{\small
    \put(-157,35){$r_1$}
    \put(-69, 105){$r_2$}
    \put(-134,12){$l_2$}
    \put(-133,91){$l_1$}
    \put(-110,38){$l_3$}
    \put(1,39){$l_4$}
    \put(-178,52){
\tikzstyle{mybox} = [draw=blue,rectangle, rounded corners, inner sep=2pt, inner ysep=2pt]
\begin{tikzpicture}
\node [mybox] (box){\color{blue} 3 4};
\end{tikzpicture}}
    \put(-32,95){
\tikzstyle{mybox} = [draw=blue,rectangle, rounded corners, inner sep=2pt, inner ysep=2pt]
\begin{tikzpicture}
\node [mybox] (box){$\color{blue}-2$};
\end{tikzpicture}}
}
\caption{The relative-slice formulation of the embedding problem in $D^4$ for two copies of the manifold
$A$ attached to the Hopf link: the curves $l_1,\ldots, l_4$ need to be sliced in the
handlebody $D^4\cup_{r_1, r_2} 0$-framed $2$-handles. (Orientations of link components are used in 
calculations in section \ref{Hopf subsection}.)}
\label{fig:Hopfdoubled}
\end{figure}

The relative-slice problems corresponding to the statements in theorem \ref{main theorem} for the Hopf link and for the Borromean rings
are shown in figures \ref{fig:Hopfdoubled} and \ref{fig:BorDouble} respectively. The circled numbers next to the components $r_1, r_2$ in figure \ref{fig:Hopfdoubled} are the indices of the slices that go over the $2$-handle attached to the curves $r_i$, this is discussed further in section
\ref{Hopf subsection} below.

In practice the slices in a solution to the relative-slice problem will be constructed by taking band sums of the
components $l_i$ with parallel copies of the curves $r_j$. These bands correspond to index $1$ critical points 
of the slices with respect to the radial Morse function on $D^4$, and parallel copies of $r_j$ bound disjoint embedded disks in the $2$-handle
attached to $r_j$. If the resulting link $L'$ is null-homotopic (in the sense discussed in section \ref{Milnor group subsection}),
the construction of the (singular) slices is completed by capping off the components of $L'$ by disjoint  disks in $D^4$.
These disks in general will have self-intersections; indeed the approach outlined here can be used either to find an obstruction or to find a solution up to link homotopy, i.e. disjoint slices which may have self-intersections,
as in section \ref{Hopf subsection}.

The Milnor group will be used to carry out this strategy to find a solution for the Hopf link in section \ref{Hopf subsection} and to find an obstruction for the Borromean rings in section \ref{Borromean section}. As indicated in section \ref{Milnor group subsection}, the Milnor group is very well suited for calculations up to link-homotopy. 
In both problems at hand, omitting one of the attaching curves for the $2$-handles gives the unlink. Moreover, any band-sum of this unlink with parallel copies of the given dotted curves (taking place in the context of the relative-slice problem) yields a homotopically trivial link. Therefore the Milnor group of the complement of the resulting band-summed link is isomorphic to the free Milnor group. The problem then is reduced to the question of whether band-summing may be performed so that the omitted component is trivial in the free Milnor group  (so that the entire link is homotopically trivial). A specific choice of band sums shows that the answer is ``yes'' for the Hopf link, and an argument analyzing Jacobi relations in the free Milnor group proves that the answer is ``no'' for the Borromean rings.

\begin{remark} \label{standard remark}
Recall the requirement, introduced before the statement of theorem \ref{main theorem},  that the embeddings $f_i\co (M, {\gamma})\hookrightarrow (D^4, S^3)$ 
in the definition of $M$-slice are ``standard'':
each embedding $f_i$ is isotopic to the original embedding  $(M, {\gamma})\subset (D^4, S^3)$. 
This requirement is reflected in the relative-slice context by the condition that the slices bounded by the
curves $l_i$ of a given $4$-manifold do not go over the $2$-handle (attached to the dotted curve) corresponding 
to the $1$-handle of the same $4$-manifold. Specifically, in figure \ref{fig:Hopfdoubled} the slices for $l_1, l_2$ should not go over the $2$-handle attached to $r_1$, and similarly $l_3, l_4$ should not go over $r_2$. (Note that without this restriction there is in fact a rather straightforward solution to this relative-slice problem.) 
In fact, the obstruction for the Borromean rings in section \ref{Borromean section} uses only a weaker consequence of this condition, 
that the curves $l_i$ {\em homologically} do not go over the $2$-handle corresponding 
to the $1$-handle of the same $4$-manifold (see Condition \ref{homological condition}). One could use this to define 
a {\em homologically standard} requirement which interpolates between arbitrary embeddings and standard embeddings, and the results of this paper
hold in this setting as well. (This may also serve as a bridge with the classical subject of slice links since  the usual slice disks are of course ``homologically standard''.)
This weaker homological condition on embeddings is not pursued further in the present paper since it does not have an immediate application in the AB slice problem which 
is the main motivation for the ``standard'' requirement.    
\end{remark}

\subsection{Proof of theorem \ref{main theorem} for the Hopf link.} \label{Hopf subsection}

\begin{figure}[t]
\includegraphics[width=5.2cm]{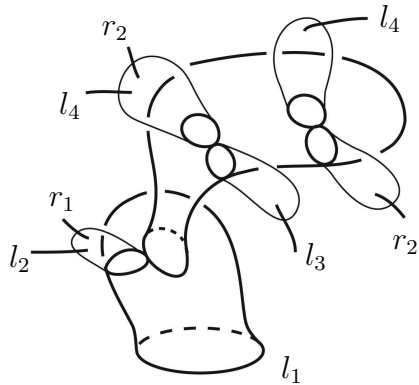}
    \put(-52,0){$l_1$}
    \put(-155,42){$l_2$}
    \put(-140,65){$r_1$}
    \put(-136,100){$l_4$}
    \put(-120,130){$r_2$}
    \put(-15,133){$l_4$}
    \put(-10,49){$r_2$}
    \put(-43,43){$l_3$}
\caption{A $2$-stage capped grope bounded by $l_1$ in figure \ref{fig:Hopfdoubled}.}
\label{fig:grope}
\end{figure}

Consider  the relative slice set-up in figure \ref{fig:Hopfdoubled}, where the dotted curves are defined in figure \ref{fig:Bnew}.
%The link $L$ formed by the components $l_1,\ldots, l_4$ is almost trivial (deleting any component yields a null-homotopic link),
%and the same statement is also true for any link formed from $L$ by taking band sums with parallel copies of $r_1, r_2$. 
To read off the word represented by $l_1$ in 
$M{\pi}_1(S^3\smallsetminus (l_2\cup l_3\cup l_4\cup r_1\cup r_2))$ it is useful to consider a $2$-stage capped grope \cite{FQ} shown
in figure \ref{fig:grope}, bounded by $l_1$: the two surface stages are embedded in the link complement while the caps
intersect the link as shown in figure \ref{fig:grope}. (There are several versions of {\em gropes} considered in the literature. Throughout this section the term ``grope'' refers to {\em half-gropes}, see definition 2.4 and figure 2.1 in \cite{FL}.)

\begin{figure}[h]
%\centering
\vspace{.2cm}
\includegraphics[height=4.5cm]{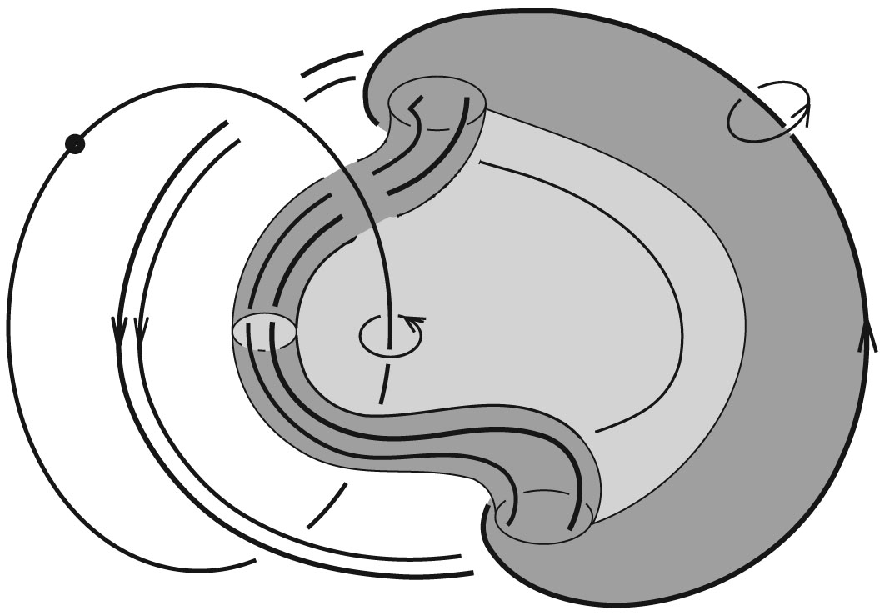} \hspace{1.5cm} 
\includegraphics[height=4.5cm]{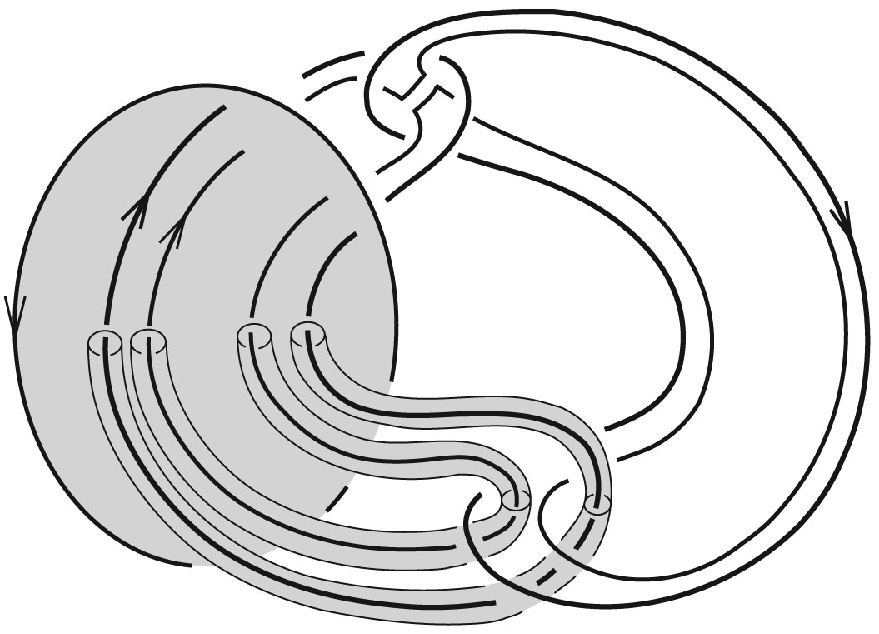}
{\small
    \put(-238,90){$l_1$}
    \put(-324,58){$\wedge$}
\put(-402,60){$l_2$}
\put(-382,60){$r_1$}
    \put(-121,120){$l_3$}
    \put(-35,60){$r_2$}
    \put(-28,110){$l_4$}
	\put(-191,60){$\wedge$}
}
\caption{A more detailed construction of the grope in figure \ref{fig:grope}. For clarity of illustration the part of the curve $r_1$ located behind the surface in the figure on the left is not shown (compare with figure \ref{fig:Bnew}).}
\label{fig:gropelocation}
\end{figure}

A more detailed construction of this grope is shown in figure \ref{fig:gropelocation}.
Specifically, the link in figure \ref{fig:Hopfdoubled} is a {\em composition} (in the sense of \cite[Theorem 2.3]{FL}) of the two links shown in figure \ref{fig:gropelocation}. Considering the standard genus one Heegaard decomposition of $S^3$, the components $l_1, l_2, r_1$ may be thought of as being contained in one of the solid tori, and the components $l_3, l_4, r_2$ are in the other solid torus of the decomposition. The first solid torus is pictured on the left in figure \ref{fig:gropelocation} as the complement of the dotted curve. 
The (genus one) first stage surface of the grope is shown in that figure. One cap for that surface intersects the components $l_2, r_1$. There is another cap, intersecting the dotted curve, visible in the picture. The second stage is obtained by removing a small meridional disk from this cap. The resulting boundary component is then filled in by the genus two surface bounded by the curve $\Lambda$ in the complement of the link in figure \ref{fig:gropelocation} on the right.

The component $l_1$ is seen to  be represented by the word 
\begin{equation} \label{l1 word}
l_1=[[m_3,m_4\cdot b]\cdot [b, m_4],m_2\cdot a]
\end{equation}
where $m_i, a,b$ denote meridians to $l_i$, $r_1$, $r_2$ respectively. 
Recall that meridians are small circles linking the components, connected by arcs to the basepoint.
The meridians, viewed as elements of the fundamental group of the link, depend on the choice of these arcs, as well as on the choice of the orientations of the small circles.
It will be clear from the proof below that the choice of connecting arcs (and 
the choice of bands) is not going to be important for the argument since the difference between various choices is measured by higher order commutators
which are trivial in the relevant Milnor group.

To fix the ambiguity with orientations, consider the orientations of the link components specified in figure \ref{fig:gropelocation} (also see figure \ref{fig:Hopfdoubled}). Orient the meridians using the ``right-hand rule'', so the linking number of each component with its meridian is $+1$, as illustrated for $l_1$ in figure \ref{fig:gropelocation}. A direct calculation shows that the exponents of all meridians in (\ref{l1 word}) are $+1$. It is worth mentioning that the techniques developed for the proof
of theorem \ref{main theorem} in this paper, both for the Hopf link and for the Borromean rings, are quite robust. They work for a family of examples generalizing 
the main example in figure \ref{fig:Bnew}. For instance, the proof goes through if the band visible in the definition of the dotted curve in figure  \ref{fig:Bnew} were twisted. In this case the exponents of the two meridians labeled $b$ in (\ref{l1 word}) would have been opposite.

Consider band sums indicated
in figure \ref{fig:Hopfdoubled}: the circled numbers are the indices of the slices (bounded by the curves $l_i$)
that go over the $2$-handles attached to $r_1$, $r_2$. Specifically, add one parallel copy of $r_1$ to both 
$l_3$ and $l_4$, and add a parallel copy of $r_2$ with reversed orientation to $l_2$. Denote the components formed by the band summing by $l'_1, l'_2, l'_3, l'_4$.
($l'_1$ equals $l_1$ since no band-summing is performed on this component.)

\begin{proposition}  \label{234 homotopically trivial}
\sl The link $(l_2', l_3', l_4')$ is homotopically trivial. 
\end{proposition}

\begin{figure}[h]
%\centering
\includegraphics[height=3.4cm]{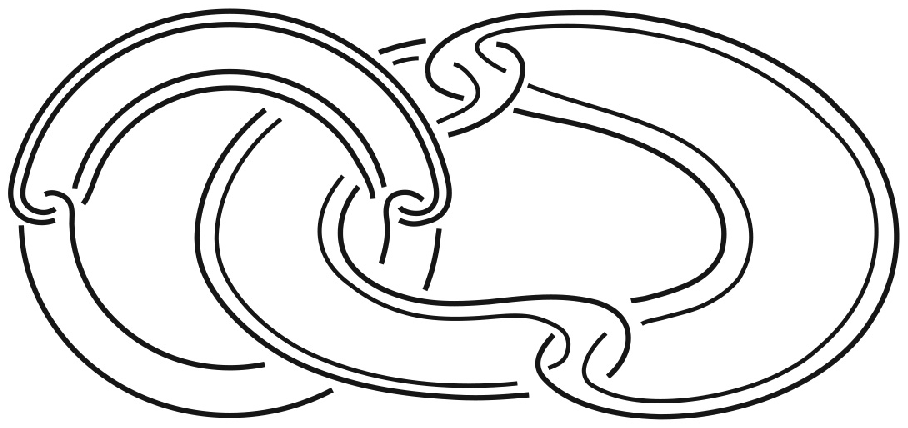}
{\scriptsize
    \put(-203,74){$\overline 3$}
    \put(-188,70){$\overline 4$}
    \put(-115,12){$\overline 2$}
    \put(-36,15){$\overline{\overline 2}$}
    \put(-106,-2){$l_3$}
    \put(-30,0){$l_4$}
	\put(-183,0){$l_2$}}
\caption{}
\label{fig:nullhomotopic}
\end{figure}

{\em Proof.} A useful tool is the {\em Half-grope lemma} \cite[Theorem 2.5]{FL} (also see 
%\cite[Theorem 2]{KT} and 
\cite[Theorem 2]{K6} for a streamlined proof). 
It states that if the components of an $n$-component link bound disjoint maps of $(n-1)$-stage gropes in $D^4$ then the link is homotopically trivial. Therefore in our context it suffices to find disjoint maps of three $2$-stage gropes into $D^4$, bounded by 
$l_2', l_3', l_4'$. Consider two parallel copies of $r_1$ in figure \ref{fig:Hopfdoubled} (without band-summing them with $l_3, l_4$).  Recall the detailed drawing of the handle diagram in the solid torus in figure \ref{fig:Bnew}. Since $l_1$ is missing from the link currently under consideration, observe that these two parallel copies of $r_1$ are isotopic within the solid torus  to the curves labeled $\overline 3, \overline 4$ in figure \ref{fig:nullhomotopic}. 

Consider the 7-component link in figure \ref{fig:nullhomotopic}.
The curve $l_2$ bounds an embedded $2$-stage grope, which can be easily located in the picture, in the complement of the other components in $S^3$. Extend these other 6 components by a product in the collar $S^3\times[0,\epsilon]$ in $D^4$. Since $l_2$ is not present in $S^3\times\{ \epsilon\}$, the 6 components form the unlink and so bound disjoint disks in $S^3\times\{ \epsilon\}$. Therefore the $7$ component link in figure \ref{fig:nullhomotopic} bounds disjoint $2$-stage gropes in $D^4$.
(By definition the disk is an $n$-stage grope, for any $n$.)
The link $(l'_2, l'_3, l'_4)$ is formed from the link in figure \ref{fig:nullhomotopic} by band-summing $l_3$ with the curve labeled $\overline 3$, $l_4$ with $\overline 4$, and also $l_2$ with $\overline 2$ and $\overline{\overline 2}$. Taking a boundary-connected sum of the $2$-stage gropes constructed above along the bands defining the band sums yields three disjoint $2$-stage gropes in $D^4$, bounded by $l'_2, l'_3, l'_4$. The half-grope lemma completes the proof of proposition \ref{234 homotopically trivial}. \qed

Recall the commutator identities (cf. \cite[Theorem 5.1]{MKS})
\begin{equation} \label{product identity}
[x,yz]\, =\, [x,z]\; [x,y]^z, \; \;\,[xz,y]\, =\, [x,y]^z\; [z,y].
\end{equation}

We will now take a brief digression to discuss a basic fact, important for the argument here and in section \ref{Borromean section}, that conjugation in the commutator identities (\ref{product identity}) and (\ref{Hall-Witt equation}) below does not affect calculations in the Milnor group in our setting. The same comment applies to conjugation that results from different choices of basepoints. The key point is that the component $l_1$ is in the third term of the lower central series of ${\pi}_1(S^3\smallsetminus(l_2'\cup l_3'\cup l_4'))$. Geometrically this is reflected in the fact that
it bounds a two stage grope in figure \ref{fig:grope}; this may be seen algebraically using the expression (\ref{l1 word}) and the identities (\ref{product identity}). 
Recall from section \ref{Milnor group subsection} that the Magnus expansion $M$ of the free Milnor group $M{\pi}_1(S^3\smallsetminus (l'_2\cup l'_3\cup l'_4))$ into the ring $R_{x_2, x_3, x_4}$ is well-defined and injective. The Magnus expansion takes any element of the third term of the lower central series to a polynomial of the form $1+$(some linear combination of monomials of length $3$ in non-repeating variables $x_2,x_3,x_4$). 
The effect of conjugation on the Magnus expansion is the addition of higher order monomials. However since the monomials in the expansion $M(l_1)$ are already of maximal length in $R_{x_2, x_3, x_4}$, any type of conjugation mentioned above does not change  $M(l_1)$. Since $M$ is injective, the element represented by $l_1$ in $M{\pi}_1(S^3\smallsetminus (l'_2\cup l'_3\cup l'_4))$ is also unchanged by conjugation. An alternative argument for this fact, not using the Magnus expansion, may be given by directly using the defining Milnor relation (\ref{eq:Milnor group}).
%The fact that conjugation does not change the relevant element in the Milnor group will apply in section \ref{Borromean section} as well, where the component $l_1$ is also in the maximal non-trivial term of the lower central series of the Milnor group of the complement of the rest of the link.

The band summing (indicated in figure \ref{fig:Hopfdoubled}) defining the link $(l'_1, l'_2, l'_3, l'_4)$  results in substitutions $a=m_3m_4$ and $b=m_2^{-1}$ in (\ref{l1 word}). Disregarding conjugation in (\ref{product identity}) and collecting commutators with distinct indices, one has
\begin{equation} \label{l1 eq}
l_1\, =\, [[m_3 ,m_4 \cdot m^{-1}_2]\cdot [m^{-1}_2,m_4], m_2\cdot m_3\cdot m_4]=\end{equation}
$$[[m_3,m_4],m_2]\cdot [[m_3,m^{-1}_2],m_4]\cdot [[m_2^{-1},m_4],m_3].$$

Using the identity $[x^{-1},y]=[y,x]^{x^{-1}}$ (where conjugation is again irrelevant), $l_1$ equals
$$[[m_3,m_4],m_2]\cdot [[m_2,m_3],m_4]\cdot [[m_4,m_2],m_3].$$

Omitting conjugation in the Hall-Witt identity \cite[Theorem 5.1]{MKS}
\begin{equation} \label{Hall-Witt equation}
[[x,y],z^x]\cdot [[z,x],y^z]\cdot [[y,z],x^y]=1,
\end{equation}

one gets the Jacobi relation 
\begin{equation} \label{Jacobi relation}
[[x,y],z]\cdot [[z,x],y]\cdot [[y,z],x]=1,
\end{equation} 

establishing that $l_1=1\in M Free_{m_2,m_3,m_4}$.

Recall from proposition \ref{234 homotopically trivial} that the link $(l'_2, l'_3, l'_4)$ is homotopically trivial. The point of the calculation above is that $l_1=1\in M{\pi}_1(S^3\smallsetminus (l'_2\cup l'_3\cup l'_4))\cong M Free_{m_2,m_3,m_4}$.
Therefore it follows from \cite[Theorem 3]{M} that $L'=(l_1, l'_2, l'_3, l'_4)$ is homotopically trivial and so its components bound disjoint maps of disks into $D^4$.
This gives a solution up to link-homotopy to the relative-slice problem for the Hopf link
which is ``standard'' in the sense of remark \ref{standard remark}.

\begin{remark} To illustrate the subtlety of the problem for the Hopf link analyzed above, 
it is interesting to note that while a link-homotopy solution is shown to exist for the relative-slice problem in figure \ref{fig:Hopfdoubled}, there are no {\em embedded} slices for the components $l_1,\ldots, l_4$ in this problem. This may be proved by finding an obstruction similar to that in \cite{K5} for a link obtained by handle slides on the link  $(l_1,\ldots, l_4)$. The condition of being relatively slice is preserved by handle slides but the condition of being relatively slice up to link homotopy is not preserved  in general.
\end{remark}

\subsection{Completion of the construction of ${\mathbf M}$.} \label{completion}
The detail that is missing from the conclusion of the proof of theorem \ref{main theorem} for the Hopf link is that the two copies
of the manifold $M$ have to be {\em embedded}, while the outcome of the argument so far is a map of two copies $A_1, A_2$ of $A$ into $D^4$
where individual $2$-handles have self-plumbings. (Moreover, as discussed in the introduction the embeddings are required to be isotopic to the original embedding
$M\subset D^4$.)
The null-homotopies produced as a result of Milnor group calculations may be realized as a sequence of standard ``elementary homotopies''
of link components.
The manifold $M$ will be defined to be $A$ with a number of self-plumbings of its two $2$-handles, determined by those of both $A_1$ and $A_2$. 

\begin{figure}[ht]
\includegraphics[height=3.5cm]{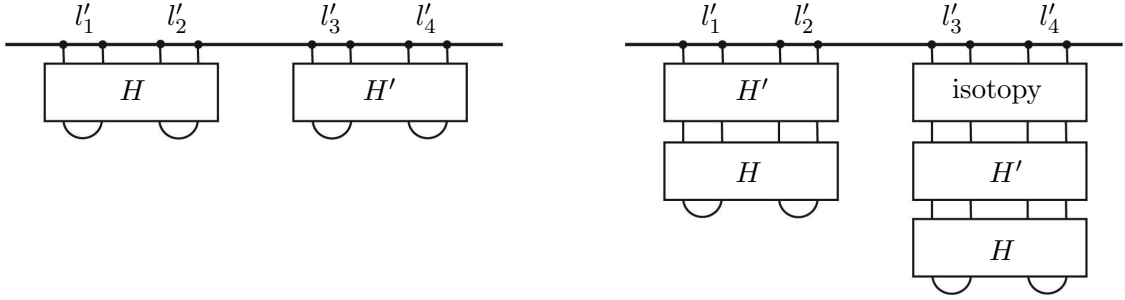}
{\small
    \put(-400,103){$l'_1$}
    \put(-365,103){$l'_2$}
    \put(-307,103){$l'_3$}
    \put(-271,103){$l'_4$}
    \put(-383,74.5){$H$}
    \put(-291,74.5){$H'$}
    \put(-163,103){$l'_1$}
    \put(-128,103){$l'_2$}
    \put(-72,103){$l'_3$}
    \put(-35,103){$l'_4$}
    \put(-150,43.5){$H$}
    \put(-150,74.5){$H'$}
    \put(-68,75.5){${\rm isotopy}$}
    \put(-54,43.5){$H'$}
    \put(-54,14){$H$}
}
\caption{A schematic picture of a preliminary link-homotopy of $L$ (left) and of a modified link homotopy (right) yielding 
two disjoint standard embeddings of $M$. The horizontal direction depicts $S^3$ and the vertical direction corresponds to the radial coordinate in $D^4$. The link components $l'_i$ are pictured as $0$-spheres in $S^3=\partial D^4$. Each box $H, H'$ denotes two disjoint level-preserving annuli, possibly with self-intersections. The top and the bottom of each of $H, H'$ individually is a two-component unlink.}
\label{homotopy fig}
\end{figure}

The precise details of the construction of $M$ are as follows. It is worth mentioning right away that the {\em standard} condition is
imposed on each individual embedding of $M$. That is, after $M$ is constructed two copies of it will be disjointly embedded with their attaching circles corresponding to the two components of the Hopf link, and the embedding of each copy will be shown to be standard after the other component is disregarded. 
Figure \ref{homotopy fig} is a concise illustration of the construction.

Each of the two links 
$(l_1, l_2, r_2)$ and $(l_3, l_4, r_1)$ is a three-component unlink, and band sums in the proof above may be easily found so that both $(l'_1, l'_2)$, $(l'_3, l'_4)$ are
two-component unlinks. Taking a band sum of $l_2$ with $r_2$ and capping off with the core of the $2$-handle attached to $r_2$ amounts to a $(1,2)$-pair of critical points of the slice for $l_2$ with respect to the radial Morse function on $D^4$. 
(This Morse function is considered on the original $4$-ball which contains the $2$-handles attached along $r_1, r_2$ and into which the manifolds $A_1, A_2$ are embedded, see the second paragraph of section \ref{relative slice subsection}. The slightly smaller $4$-ball where the relative-slice problem is being considered is obtained by removing the collars on the attaching regions of $A_1, A_2$.)
Since 
$(l_1, l_2, r_2)$ is the unlink, this pair of critical points can be canceled and (ignoring the components $l_3, l_4$) the result is an isotopy of $(l_1,l_2)$ in $S^3\times I$. (Compare with figure 15 in \cite{K4}.)
Similarly, (considering just the last two components) band summing $l_3$ and $l_4$ with copies of $r_1$ may be realized instead as an isotopy of $l_3, l_4$. 

We summarize the set-up: the $4$-component link $L':=(l'_1,\ldots,l'_4)$ is null-homotopic, and the two component sub-links $(l'_1, l'_2), (l'_3, l'_4)$ are 
individually unlinks. We will next construct specific null-homotopies of the components of $L'$.  Start with any link homotopy from $L'$ to the $4$-component unlink.
Rather than capping them off by disks right away, let the first two components move by an isotopy away from the last two components.
For reasons which will be clear below, next we run the entire link-homotopy of $l'_3, l'_4$ backwards, while the first two components stay fixed as unknots away from $l'_3, l'_4$.
The result is a $4$-component unlink which may be capped off by disks. Denote the resulting link homotopies (disjoint annuli with self-intersections) of
$l'_1, l'_2$ by $H$ and $l'_3, l'_4$ by $H'$, see the diagram on the left in figure \ref{homotopy fig}. To prepare for the following step of the construction, note that 
both $H, H'$ are level-preserving singular maps $(S^1\sqcup S^1)\times [0,1]\longrightarrow S^3\times [0,1]$ giving two-component unlinks at times $0,1$. Their singularities consist of finitely many double point self-intersections, and they are isotopies when restricted to sub-intervals of $[0,1]$ not containing singular points. Although $H'$ is defined using $l'_3, l'_4$, it can be ``applied'' to any two-component unlink in $S^3\times\{ 0\}$, for example to $l'_1, l'_2$.

%There exist link-homotopies of the components 
%$l'_2, l'_3, l'_4$ so that $l'_1$ is null-homotopic in the complement of the resulting link in $S^3$. (For example,  they may be chosen to be a sequence of
%standard ``finger move'' homotopies which are geometric realizations of the defining relations (\ref{eq:Milnor group}) of the Milnor group $M{\pi}_1(S^3\smallsetminus (l'_2, l'_3, l'_4))\cong M{\rm Free}_{m_2,m_3,m_4}$.) Rather than capping $l'_1$ by a singular disk in the complement
%of the other components, let it move by a link-homotopy to an unknot away from the other $3$-components. 
%The link formed by the components $l'_2, l'_3, l'_4$ at this time is a result of a link homotopy applied to the last three components of $L'$, so it is null-homotopic. Apply the preceding argument to this three-component link: let $l'_3, l'_4$ move by a link-homotopy until $l'_2$ 
%may be homotoped in their complement to an unknot away from them. 

The result so far is insufficient for the definition of $M$ since $H, H'$ may be non-isotopic in $D^4$. 
A suitable link-homotopy of $L'$ is constructed instead as follows. As shown on the right in figure \ref{homotopy fig}, as a preliminary step start by applying the link homotopy $H'$ to the components $l'_1, l'_2$, while letting $l'_3, l'_4$ move by an isotopy
in their complement. More precisely, recall from the previous paragraph that $H'$ is a level-preserving map $(S^1\sqcup S^1)\times I\longrightarrow S^3\times I$ which is an isotopy at generic times and whose singularities consist of double point self-intersections at finitely many times. Given two components, $l'_3$ and $l'_4$, in the complement of $l'_1, l'_2 $ in $S^3\times\{ 0 \}$ there exists an isotopy moving them in the complement of $H'$: a level-preserving embedding $(S^1\sqcup S^1)\times I\hookrightarrow S^3\times I\smallsetminus \, {\rm image}(H')$. Indeed, each double point of $H'$ is described by a movie where two strands of either $l'_1$ or $l'_2$ move by an isotopy in $S^3$ and intersect in a single point. By general position, since $l'_3, l'_4$ are $1$-dimensional submanifolds of $S^3$, they may be kept disjoint from $H'$ during this movie. This completes the argument for the existence of an isotopy of $l'_3, l'_4$ in the complement of $H'$. 

Since $H'$ was constructed as a link-homotopy followed by its reverse, the outcome of the previous paragraph is an identical copy of $L'=(l'_1, \ldots, l'_4)$. Now run 
the link-homotopy $H, H'$ of $L'$. The result is a $4$-component unlink. Cap off the first two components with disks and 
apply $H$ to the $3$-rd and $4$-th components. Finally, cap off the last two components. Now the null-homotopy of $(l'_1, l'_2)$ is isotopic to the null-homotopy of $(l'_3, l'_4)$: up to isotopy each one is $H'$ followed by $H$ and then capped off with standard disks, figure \ref{homotopy fig}.
Define $M$ to be either of the two embeddings  of the manifold $A$ in figure \ref{fig:Bnew} with self-plumbings of the $2$-handles
constructed above.
This yields an embedding of two standard copies of $M$ as required in the part of theorem \ref{main theorem} concerning the Hopf link.

\section{An obstruction for the Borromean rings} \label{Borromean section}

This section completes the proof of theorem \ref{main theorem} by showing 
that the Borromean rings do not bound disjoint standard embeddings of three copies of $(M,{\gamma})$ in $D^4$.
This argument is important in the context of the A-B slice problem, discussed in section \ref{AB slice section}.

\begin{figure}[ht]
%\centering
\vspace{.2cm}
\includegraphics[width=6cm]{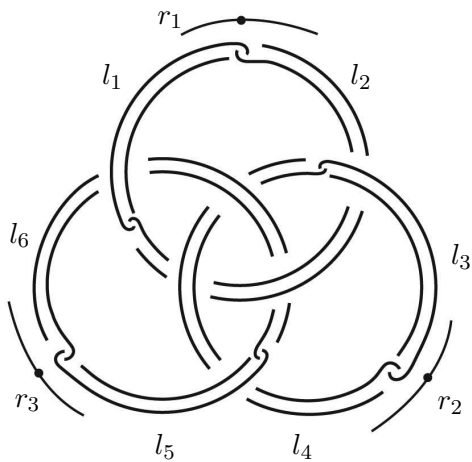}
{\small
    \put(-170,73){$l_6$}
    \put(-135,133){$l_1$}
    \put(-41,133){$l_2$}
    \put(-115,-7){$l_5$}
    \put(-63,-7){$l_4$}
    \put(-3,65){$l_3$}
    \put(-168,11){$r_3$}
    \put(-8,10){$r_2$}
    \put(-114,156){$r_1$}}
    \vspace{.45cm} \caption{The relative-slice problem for embedding three copies of $A$ on the Borromean rings: do the components
$l_1,\ldots, l_6$ bound disjoint disks in $D^4\cup_{r_1,r_2,r_3}2{\rm-handles}$? The dotted curves $r_1, r_2, r_3$ are shown only schematically,
as explained in figure \ref{fig:Bnew}.}
\label{fig:BorDouble}
\end{figure}

The relative slice formulation  corresponding to the embedding
problem for three copies of the manifold $A$ in section \ref{construction section} is shown in figure \ref{fig:BorDouble}.
(The manifold $M$ in the statement of theorem \ref{main theorem} was defined in section \ref{completion}
as $A$ with a number of self-plumbings of its $2$-handles. As noted in section \ref{Milnor group subsection} 
these self-plumbings do not affect the Milnor group
arguments given below.)

Suppose to the contrary that the Borromean rings bound disjoint standard embeddings of three copies $A_i$, $i=1,2,3$ of $A$; equivalently
assume the link in figure \ref{fig:BorDouble} is relatively slice, subject to the ``standard'' condition discussed
in remark \ref{standard remark}. We specify a weak consequence of this condition, sufficient for the proof of theorem \ref{main theorem} for the Borromean rings:

\begin{condition} \label{homological condition}
The slices for $l_1,l_2$ do not homologically  go over the $2$-handle attached to $r_1$, and similarly $l_3, l_4$  do not  go over $r_2$, and $l_5, l_6$ do not go over $r_3$. 
Here, a slice {\em does not homologically go over a $2$-handle} if its (${\mathbb Z}$-valued) algebraic intersection number with the co-core of the $2$-handle is zero.
\end{condition}

It will be shown that the components $l_i$ do not bound disjoint singular disks subject to this restriction.
(Without this restriction, it is not hard to find a solution to this relative-slice problem.)
This condition is important in the A-B slice problem
\cite{FL}, see section \ref{AB slice section}.

%Denote by $D_i$ the slice bounded by $l_i$, $i=1,\ldots, 6$, and let $D$ be the union of the first five slices,
%$D=\cup_{i=1,\ldots, 5} D_i$. Consider $$X:= (D^4\cup_{r_1,r_2,r_3}2{\rm-handles})\smallsetminus D.$$
%Denote meridians to the components $l_i$ by $m_i$ 
%The first homology $H_1(X)$ is generated by $m_1,\ldots, m_5$; abusing the notation we will 
%view $\{m_i\}$ as based loops in $X$ normally generating ${\pi}_1(X)$ and also as generators of $M{\pi}_1(X)$.

\begin{figure}[ht]
%\centering
\includegraphics[height=3.5cm]{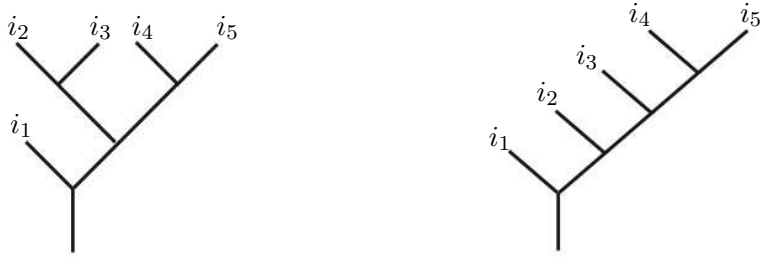}
{\small
    \put(-288,53){$i_1$}
    \put(-289,90){$i_2$}
    \put(-258,90){$i_3$}
    \put(-242,90){$i_4$}
    \put(-210,90){$i_5$}
    \put(-107,48){$i_1$}
    \put(-89,65){$i_2$}
    \put(-74,79){$i_3$}
    \put(-54,95){$i_4$}
    \put(-12,95){$i_5$}}
\caption{The tree on the left encodes the bracketing in the commutators of the form $[m_{i_1},[[m_{i_2},m_{i_3}],[m_{i_4},m_{i_5}]]]$ arising
in the link in figure \ref{fig:BorDouble}, and the tree on the right encodes commutators $[m_{i_1},[m_{i_2},[m_{i_3},[m_{i_4},m_{i_5}]]]]$ spanning the vector space $V$.}
\label{fig:trees}
\end{figure}

Consider the 
$120$-dimensional vector space $V$ over ${\mathbb Q}$, formally spanned by all commutators of the form 
$[m_{i_1},[m_{i_2},[m_{i_3},[m_{i_4},m_{i_5}]]]]$ in five non-repeating variables,
so the indices $({i_1},{i_2},{i_3},{i_4},{i_5})$ range over all permutations of $\{2,\ldots,6\}$. Omitting the curves $r_1, r_2, r_3$ in figure \ref{fig:BorDouble}, 
the remaining link is the Bing double of the Borromean rings. Using the orientations in figure  \ref{fig:BorDouble} and the corresponding choice of meridians
(discussed in section \ref{Hopf subsection}), one checks that the component $l_1$ represents the commutator 
\begin{equation} \label{l6} 
l_1=[m_2,[[m_3,m_4],[m_5,m_6]]].
\end{equation}

in the complement of the other five components $l_2, \ldots, l_6$.
An observation important for the proof below is that any band sum of the curves $\{l_i\}$ and parallel copies of the dotted curves $r_j$ in the link in figure  
\ref{fig:BorDouble}
gives rise to (products of) commutators of the form
\begin{equation} \label{CommutatorForm}
[m_{i_1},[[m_{i_2},m_{i_3}],[m_{i_4},m_{i_5}]]]
\end{equation}
where $i_1,\ldots, i_5$ is a permutation of
$\{2,\ldots, 6\}$,  see figure \ref{fig:trees}.  Any such commutator is equivalent under anti-symmetry relation applied to $[[m_{i_2},m_{i_3}],[m_{i_4},m_{i_5}]]$ to one of 15 that appear in the statement
of lemma \ref{intersection prop} below.
Let $U$ be the subspace of $V$ generated by these 15 commutators. Denote by $J$
the subspace of $V$ spanned by the Jacobi and anti-commutation relations.

\begin{remark} 
Let $F$ denote the free group $Free_{m_2,\ldots,m_6}$.
The quotient $V/J$ is isomorphic to $(MF)^5\otimes{\mathbb Q}$, where $(MF)^5$ denotes the $5$-th term of the lower central series of the Milnor group $MF$.
\end{remark}

\begin{remark} 
The notation will alternate between the {\em product} of commutators, considered in the Milnor group $MF$ and the {\em sum} of commutators considered in the abelian group $(MF)^5$. (Note that $(MF)^6$ is trivial.) This should not cause any confusion.
\end{remark}

\begin{lemma} \label{intersection prop} \sl
The intersection $U\cap J$ is a $1$-dimensional subspace of $V$ spanned by the product of 15 commutators 
$w:=$
$$[m_2,[[m_3,m_4],[m_5,m_6]]]
\cdot 
[m_2,[[m_5,m_3],[m_4,m_6]]] 
\cdot 
[m_2,[[m_4,m_5], [m_3,m_6]]] $$
$$\cdot
[m_3,[[m_4,m_2],[m_5,m_6]]]
\cdot 
[m_3,[[m_2,m_5],[m_4,m_6]]]
 \cdot 
[m_3,[[m_5,m_4], [m_2,m_6]]]   $$
$$\cdot
[m_4,[[m_2,m_3],[m_5,m_6]]]
\cdot 
[m_4,[[m_5,m_2],[m_3,m_6]]] 
\cdot 
[m_4,[[m_3,m_5],[m_2,m_6]]]$$
$$\cdot
[m_5,[[m_3,m_2],[m_4,m_6]]]
\cdot 
[m_5,[[m_2,m_4],[m_3,m_6]]] 
\cdot 
[m_5,[[m_3,m_4], [m_2,m_6]]]$$
$$
\cdot[m_6,[[m_2,m_3],[m_4,m_5]]]
\cdot 
[m_6,[[m_4,m_2],[m_3,m_5]]] 
\cdot 
[m_6,[[m_2,m_5],[m_3,m_4]]]$$

Therefore dim$(U/(J\cap U))=14$.
\end{lemma}

The element $w$ is certainly in $U$; an explicit calculation in section \ref{commutator list} using the Jacobi relations (\ref{Jacobi relation}) 
shows that $w$ is also in $J$.

The following proposition describes a convenient basis of the space of commutators
in non-repeating variables modulo the Jacobi and anti-symmetry relations. (A directly analogous statement works for $n$-fold commutators for arbitrary $n$; the result is stated below in the case $n=5$ relevant for the current proof.) Of course the index ``6'' in the statement below can be replaced
by any fixed index in $\{ 2,\ldots, 6\}$. The proof of lemma \ref{intersection prop} will follow from the following proposition. 

\begin{proposition} \label{basis} \sl
The collection of commutators of the form $[m_{i_1}, [m_{i_2},[m_{i_3},[m_{i_4},m_6]]]]$, where the right-most index is $6$ and $i_1,\ldots i_4$ range over all permutations of $2,\ldots,5$,
forms a basis of $V/J$.
\end{proposition}

{\em Proof.} To show that the commutators in the statement span $V/J$, start with any commutator $[m_{i_1}, [m_{i_2},[m_{i_3},[m_{i_4},m_{i_5}]]]]$ and use the anti-symmetry relation to shift $m_5$ to the right-most position. In general this will change the bracketing pattern of the commutator. Now the Jacobi relations (\ref{Jacobi relation}) may be used to rebracket and get a product of commutators
of the form 
\begin{equation} \label{basic commutators}
[m_{i_1}, [m_{i_2},[m_{i_3},[m_{i_4},m_6]]]].
\end{equation}

To show that these elements are linearly independent in $V/J$, consider the Magnus expansion (\ref{MagnusMilnor}) and note that  the monomial $x_{i_1} x_{i_2}x_{i_3} x_{i_4} x_6$ is present only in the Magnus expansion of the commutator $[m_{i_1}, [m_{i_2},[m_{i_3},[m_{i_4},m_6]]]]$ . \qed

\medskip

{\em Proof of lemma \ref{intersection prop}.} It follows from proposition \ref{basis} that $V/J$ is 24-dimensional. 
Recall that $U$ is a 15-dimensional subspace of $V$ spanned by the commutators that appear in the statement of 
lemma \ref{intersection prop}. Using the Jacobi and anti-symmetry relations, 12 among these 15 commutators are seen to be products
of two ``basic'' commutators (\ref{basic commutators}). In total these $24=12\cdot 2$ basic commutators form the basis of $V/J$ and so clearly 
no linear combination of the first 12 commutators listed in lemma \ref{intersection prop} may intersect $J$ non-trivially. There are also 3 ``complicated''
commutators which appear last in lemma \ref{intersection prop}. Each of them is a product of 8 commutators of the form (\ref{basic commutators}), see section \ref{commutator list}, and the total $24=8\cdot 3$ of them again are a basis of $V/J$. There is only one non-trivial linear combination among the
15 commutators, the one stated in the lemma. \qed

To conclude the proof that the link in figure \ref{fig:BorDouble} is not relatively slice,
recall from (\ref{l6}) that without the components $\{r_j\}$ the curve $l_1$ reads off the first commutator in the definition of $w$ in lemma 
\ref{intersection prop}. Let $w'$ denote the rest of them, so $w=[m_2,[[m_3,m_4],[m_5,m_6]]] \cdot w'$. Then for the relative-slice problem in figure \ref{fig:BorDouble} to have a solution, the slices going over $r_1, r_2, r_3$ must ``account for'' the product $w'$. 

\begin{proposition}  \label{not in image} \sl
The element $w'$ cannot be realized by band sums of the relative-slice problem in figure \ref{fig:BorDouble}.
\end{proposition}

The proof of this proposition consists of an explicit check that the corresponding system of equations does not have a solution.
In this problem there are 14  equations of degrees 2 and 3 in 12 variables, explicitly written down in the Appendix, see section \ref{proof subsection}. 
The fact that there are no integral solutions is verified using Mathematica (section \ref{Mathematica calculation}), although using symmetries of the equations it is in fact possible to check this fact by hand. It is interesting to note that there are {\em rational} solutions which do not seem to have a geometric interpretation in the context of the relative-slice problem.

It is instructive to compare proposition \ref{not in image} with section \ref{obstruction subsection} and figure \ref{fig:BorDoublePrime} where a solution
is shown to exist for a closely related problem.

We summarize the argument proving Theorem \ref{main theorem} for the Borromean rings. Suppose there exist disjoint standard embeddings of three copies of the manifold $M$ bounded by the Borromean rings. Therefore there exists a solution to the relative slice problem in figure \ref{fig:BorDouble}, subject to Condition \ref{homological condition}. 
Then, as discussed in  section \ref{relative slice subsection}, there exist band sums of the components $l_1,\ldots, l_6$ with parallel copies of $r_1, r_2, r_3$ 
(again subject to Condition \ref{homological condition}) such that the resulting link $(l'_1, \ldots, l'_6)$ is homotopically trivial. A direct generalization of Proposition 
\ref{234 homotopically trivial} shows that for any such band sum, omitting the
first component, the rest of the link $(l'_2, \ldots, l'_6)$ is homotopically trivial. Proposition \ref{not in image}  implies that for any band sum, $l'_1$ is non-trivial in 
$M{\pi}_1(S^3\smallsetminus (l'_2, \ldots, l'_6))\cong MFree_{2,\ldots,6}$. Therefore for any band sum the link $(l'_1, \ldots, l'_6)$ is not null homotopic.
This contradiction concludes the proof of Theorem \ref{main theorem}.
\qed

\begin{remark}
An alternative to converting slices into band sums
is to use the Milnor group of the $4$-manifold $X^4:=D^4\cup_{r_1, r_2, r_3}2$-handles$\smallsetminus$(slices for $l_1, \ldots, l_6$), see \cite[Section 4]{K5}. 
The Milnor group calculations in the proof of Theorem \ref{main theorem} in both contexts are identical.
\end{remark}

\section{The A-B slice problem.} \label{AB slice section}

\subsection{Analysis of the decomposition ${\mathbf{D^4=A\cup B}}$} \label{ABdecomposition subsection}
\begin{figure}[ht]
%\centering
\includegraphics[height=3.7cm]{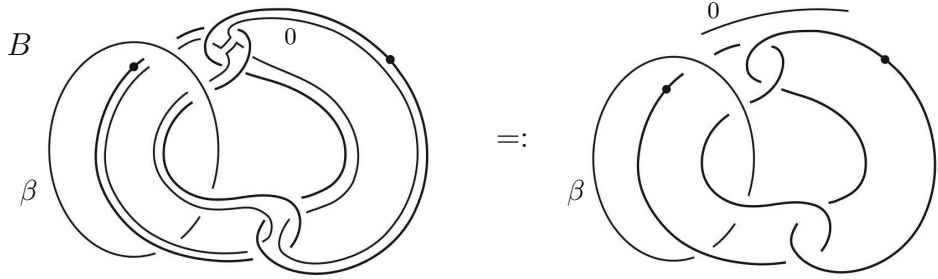}
{\small
    \put(-350,30){$\beta$}
    \put(-143,30){$\beta$}}
    \put(-355,85){$B$}
    \put(-170,50){$=:$}
{\scriptsize
    \put(-250,89){$0$}
    \put(-90,99){$0$}}
\caption{A Kirby diagram of the $4$-manifold $B=D^4\smallsetminus A$. As in figure \ref{fig:Bnew}, the picture on the right will serve as a short-hand notation for $B$.}
\label{fig:Anew}
\end{figure}

Consider the complement of the manifold $A$ in figure \ref {fig:Bnew}, $B=D^4\smallsetminus A$. 
A handle decomposition for $B$ (with the attaching curve $\beta$)  may be obtained from that of $A$ by exchanging the
$1$- and $2$-handles as explained in \cite{FL}, see figure \ref{fig:Anew}. The question relevant from the perspective of the A-B slice problem is whether the Borromean rings bound disjoint embeddings of three copies of $A$, and also disjoint embeddings of three copies of $B$ where the embeddings are all {\em standard}, see remark \ref{standard remark}. (The origin of the ``standard'' condition on embeddings is in the formulation of the A-B slice problem \cite{F2, FL}: it reflects the covering group action on the $4$-ball, which is predicted by the $4$-dimensional topological surgery conjecture.)

It was proved in section \ref{Borromean section} that there are no disjoint embeddings of three copies of $A$.
%(although a solution exists for
%a closely related problem, see section \ref{obstruction subsection} and figure \ref{fig:BorDoublePrime}!) 
We will next examine the embedding problem for three copies of $B$ with the boundary condition given by the Borromean rings. It is worth noting that the Hopf link is not $B$-slice. Observe in figure \ref{fig:Anew} that the attaching curve $\beta$ bounds a surface in $B$. Therefore the non-triviality of the linking number of the Hopf link implies that its components do not bound disjoint copies of $B$. In fact, this argument immediately generalizes to show that the Hopf link is not $A$-$B$ slice: the attaching curve of one of the two sides of any decomposition $D^4=A\cup B$ bounds rationally. Therefore the linking number is an obstruction. The interested reader will find related in-depth discussion in \cite{FK}.

\begin{figure}[ht]
\includegraphics[height=7cm]{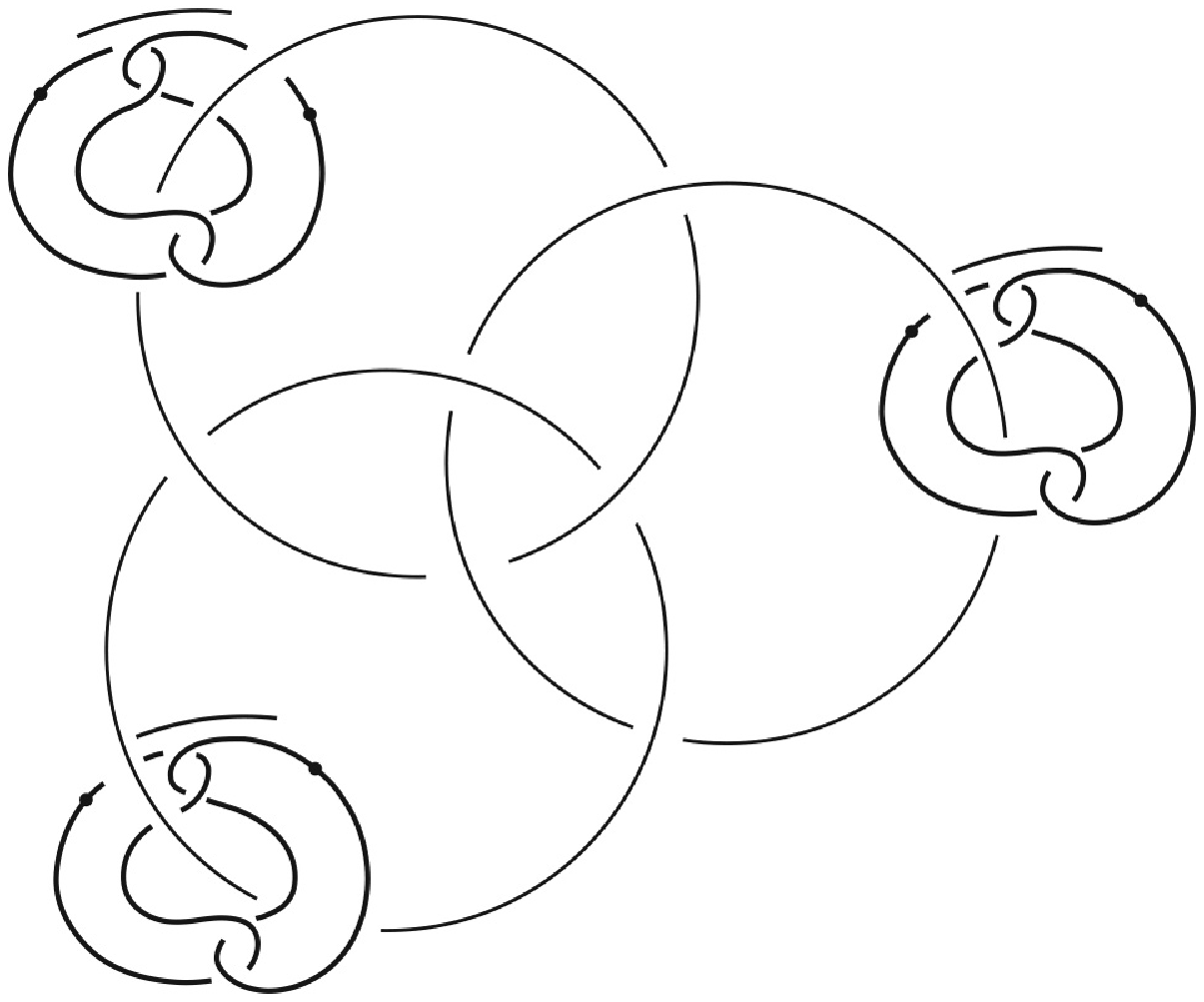}
{\small    
    \put(-126,189){$l_1$}
    \put(-60,52){$l_2$}
    \put(-132,13){$l_3$}
    \put(-217,200){$l_4$}
    \put(-36,154){$l_5$}
    \put(-202,62){$l_6$}
    \put(-180.5,168){
\tikzstyle{mybox} = [draw=blue,rectangle, rounded corners, inner sep=2pt, inner ysep=2pt]
\begin{tikzpicture}
\node [mybox] (box){\color{blue}3 5 6 };
\end{tikzpicture}}
    \put(-277,152){
\tikzstyle{mybox} = [draw=blue,rectangle, rounded corners, inner sep=2pt, inner ysep=2pt]
\begin{tikzpicture}
\node [mybox] (box){\color{blue}2 5 6};
\end{tikzpicture}}
}
{\small
	\put(-8,119){
\tikzstyle{mybox} = [draw=blue,rectangle, rounded corners, inner sep=2pt, inner ysep=2pt]
\begin{tikzpicture}
\node [mybox] (box){\color{blue}1 4 6};
\end{tikzpicture}}
	\put(-96,92){
\tikzstyle{mybox} = [draw=blue,rectangle, rounded corners, inner sep=2pt, inner ysep=2pt]
\begin{tikzpicture}
\node [mybox] (box){\color{blue}3 4 6};
\end{tikzpicture}}
	\put(-268,13){
\tikzstyle{mybox} = [draw=blue,rectangle, rounded corners, inner sep=2pt, inner ysep=2pt]
\begin{tikzpicture}
\node [mybox] (box){\color{blue}2 4 5};
\end{tikzpicture}}
	\put(-176,37){
\tikzstyle{mybox} = [draw=blue,rectangle, rounded corners, inner sep=2pt, inner ysep=2pt]
\begin{tikzpicture}
\node [mybox] (box){\color{blue}1 4 5 };
\end{tikzpicture}}}
\caption{The relative slice problem for three copies
of $B$ (figure \ref{fig:Anew}) on the Borromean rings.}
					\label{fig:ABsliceAside}
					\end{figure}

As in section \ref{Borromean section} the embedding question for the Borromean rings is reformulated as a relative-slice problem, shown in figure \ref{fig:ABsliceAside}. 
The problem is to find disjoint disks (possibly with self-intersections) for the
link components labeled $1$-$6$ in the handlebody $D^4\cup$ six zero-framed $2$-handles attached along the dotted curves. 
The sought solution has to satisfy the {\em standard} assumption discussed in remark \ref{standard remark}.
This means that the slices for $l_1, l_4$ are not allowed to go over the two $2$-handles attached to the dotted curve on the upper left, and similarly $l_2, l_5$ should not go over the two $2$-handles on the right, and $l_3, l_6$ should not go over the $2$-handles on the lower left.

The circled indices next to the dotted component in figure \ref{fig:ABsliceAside} indicate band-sums such that all $\bar\mu$-invariants of the resulting link of length $\leq 3$ vanish.
The original obstruction for the Borromean rings, $\bar{\mu}_{123}$, is of length 3, so the ``primary'' obstruction is killed.
To find a link-homotopy solution to this relative slice problem one needs all $\bar\mu$-invariants up to order 6 to be zero. 
Unlike the settings in sections \ref{Hopf section}, \ref{Borromean section}, conjugation corresponding to the choice of meridians and bands
is important here. 
It seems likely that a careful choice of bands should give rise to a link homotopy solution, but at present this is an open question.

\subsection{An obstruction for a family of examples} \label{obstruction subsection}

\begin{figure}[ht]
%\centering
\vspace{.2cm}
\includegraphics[height=3.7cm]{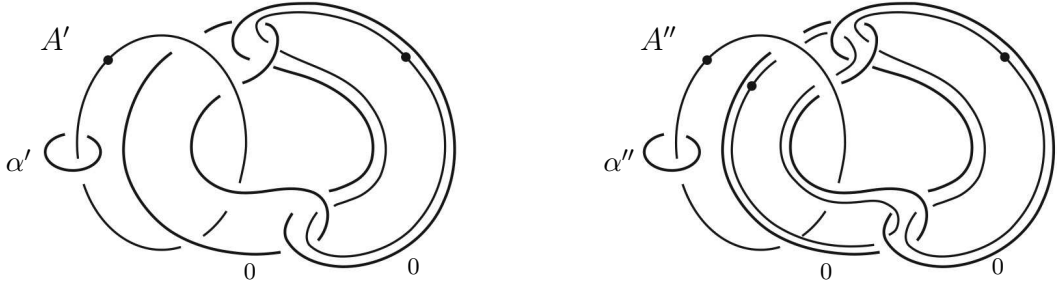}
{\small
    \put(-403,39){${\alpha}'$}
    \put(-177,39){${\alpha}''$}}
    \put(-390,85){$A'$}
    \put(-163,85){$A''$}
{\scriptsize
    \put(-313,-2){$0$}
    \put(-251,0){$0$}
    \put(-95,-2){$0$}
    \put(-30,0){$0$}}
    \vspace{.45cm} \caption{Two variations of the example from section \ref{construction section}.}
\label{fig:more examples}
\end{figure}

Consider the submanifolds $A', A''\subset D^4$ shown in figure \ref{fig:more examples}, closely related 
to the submanifold $A$ constructed in section \ref{construction section}, and consider the corresponding
decompositions $D^4=A'\cup B'$, $D^4=A''\cup B''$.

\begin{figure}[ht]
%\centering
\vspace{.2cm}
\includegraphics[width=6.1cm]{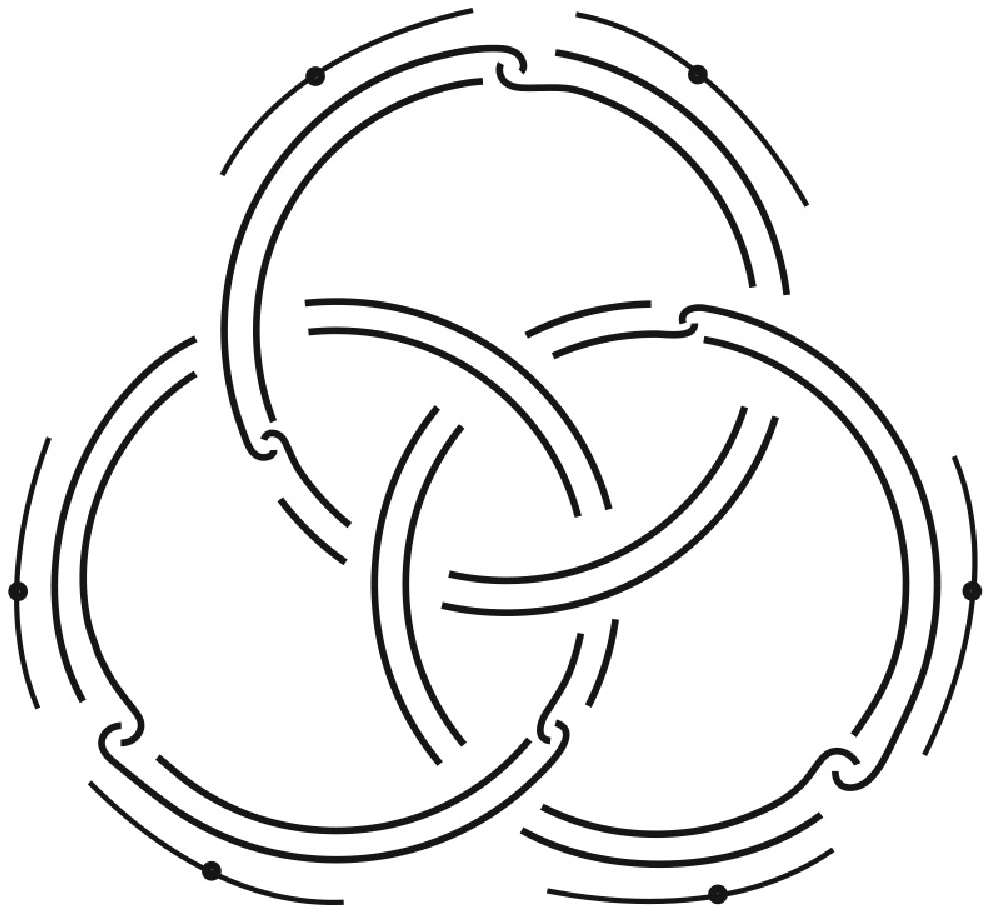}
{\small
%    \put(-199,65){$\setlength{\fboxsep}{2pt}\fbox{1\,2\,4}$}
%    \put(-168,133){$\setlength{\fboxsep}{2pt}\fbox{3\,4\,5\,6}$}
%    \put(-27,133){$\setlength{\fboxsep}{2pt}\fbox{3\,4}$}
%    \put(-125,-12){$\setlength{\fboxsep}{2pt}\fbox{1\,2\,4}$}
%    \put(-66,-12){$\setlength{\fboxsep}{2pt}\fbox{2\,5\,6}$}
%    \put(5,65){$\setlength{\fboxsep}{2pt}\fbox{2\,5\,6}$}
    \put(-209,65){
\tikzstyle{mybox} = [draw=blue,rectangle, rounded corners, inner sep=2.5pt, inner ysep=2.5pt]
\begin{tikzpicture}
\node [mybox] (box){\color{blue}1 2 4};
\end{tikzpicture}}
    \put(-179,133){
\tikzstyle{mybox} = [draw=blue,rectangle, rounded corners, inner sep=2.5pt, inner ysep=2.5pt]
\begin{tikzpicture}
\node [mybox] (box){\color{blue}3 4 5 6};
\end{tikzpicture}}
    \put(-42,133){
\tikzstyle{mybox} = [draw=blue,rectangle, rounded corners, inner sep=2.5pt, inner ysep=2.5pt]
\begin{tikzpicture}
\node [mybox] (box){\color{blue}3 4};
\end{tikzpicture}}
    \put(-145,-17){
\tikzstyle{mybox} = [draw=blue,rectangle, rounded corners, inner sep=2.5pt, inner ysep=2.5pt]
\begin{tikzpicture}
\node [mybox] (box){\color{blue}1 2 4};
\end{tikzpicture}}
    \put(-76,-17){
\tikzstyle{mybox} = [draw=blue,rectangle, rounded corners, inner sep=2.5pt, inner ysep=2.5pt]
\begin{tikzpicture}
\node [mybox] (box){\color{blue}2 5 6};
\end{tikzpicture}}
    \put(-7,65){
\tikzstyle{mybox} = [draw=blue,rectangle, rounded corners, inner sep=2.5pt, inner ysep=2.5pt]
\begin{tikzpicture}
\node [mybox] (box){\color{blue}2 5 6};
\end{tikzpicture}}
    \put(-149,66){$l_6$}
    \put(-111,124){$l_1$}
    \put(-64,125){$l_2$}
    \put(-120,20){$l_5$}
    \put(-57,20){$l_4$}
    \put(-27,60){$l_3$}}
    \vspace{.45cm} \caption{The dotted curves of three copies of $A''$ 
are shown only schematically (see the Kirby diagram on the right in figure \ref{fig:more examples}).
The circled indices show how the components $l_1, \ldots, l_6$ go over the $2$-handles attached to the dotted curves, yielding
the element $w$ from lemma \ref{intersection prop} and therefore solving the relative-slice problem up to link-homotopy.}
\label{fig:BorDoublePrime}
\end{figure}

A solution up to link-homotopy to the relative-slice problem for three copies of $A''$ on the Borromean rings is given in figure \ref{fig:BorDoublePrime}. 
(Note the assymetry of labels going
over the top two components: some assymetry is necessary since a solution does not exist for three copies of $A$, according to theorem 
\ref{main theorem}.) The reader is encouraged to check that the resulting product of commutators representing $l_6$ in non-repeating variables indeed equals 
the element $w$ from lemma \ref{intersection prop}! 
It is not difficult to check using $\bar{\mu}_{123}$ that the Borromean rings do not support standard embeddings of three copies of the other side of this decomposition, $B''$.

The decomposition $D^4=A'\cup B'$ was the subject of the papers \cite{K4,K5}. It is shown in \cite{K5} that the Borromean rings
do not bound three standard copies of $A'$.

\subsection{Summary} \label{conclusion section}
Given a decomposition of the $4$-ball, $D^4=A\cup B$, into two codimension zero submanifolds
where the ``attaching curves'' ${\alpha}, {\beta}$ of $A, B$ form the Hopf link in $S^3=\partial D^4$,
an important question in the A-B slice program is to determine whether there is necessarily an obstruction
to the Borromean rings being both $A$-slice and $B$-slice in the sense of theorem \ref{main theorem}. 
The side which carries an obstruction (if there is one) is called ``strong'', see
\cite{FL} and also \cite{K4, K5}. The goal is to determine a strong side for any 
decomposition $D^4=A\cup B$.

To summarize the results of this paper in the context of the AB slice problem, in each of the examples considered here, 
$D^4=A\cup B=A'\cup B'=A''\cup B''$ (figures \ref{fig:Bnew}, \ref{fig:more examples}), one of the two sides is found to be ``strong''.
A novel type of an obstruction is used here: the key to deciding which side is strong is whether 
the element $w$ in lemma \ref{intersection prop} is in the image of the relator curves on the relevant (in our notation, $A$)-side. 
The work presented here admits an immediate generalization giving rise to an obstruction for 
an infinite family of decompositions by further Bing doubling the links describing the Kirby diagrams in figures
\ref{fig:Bnew}, \ref{fig:more examples}.

\section{Appendix: detailed calculations.}

\subsection{A detailed commutator list for the proof of lemma \ref{intersection prop}} \label{commutator list}
To provide a verification that the element $w$ in lemma \ref{intersection prop} is in the subspace $J\subset V$, 
in other words that $w$ is trivial modulo the Jacobi and antisymmetry relations, listed below is an expression of each 
commutator in the definition of $w$ as a linear combination of the commutators forming the basis of $V/J$ in proposition \ref{basis}.
This list also makes it clear that the dimension of $U\cap J$ in the statement of lemma \ref{intersection prop} is precisely $1$ (and not greater).

\renewcommand{\arraystretch}{1.3}
\begin{tabular}{ l l}
$[m_2,[[m_3,m_4],[m_5,m_6]]] = \; $ &  $ [m_2,[m_3,[m_4,[m_5,m_6]]]] - [m_2,[m_4,[m_3,[m_5,m_6]]]]$ \\
$[m_2,[[m_5,m_3],[m_4,m_6]]] = \; $ &  $[m_2,[m_5,[m_3,[m_4,m_6]]]] - [m_2,[m_3,[m_5,[m_4,m_6]]]]$ \\
$[m_2,[[m_4,m_5], [m_3,m_6]]] = \; $ &  $[m_2,[m_4,[m_5,[m_3,m_6]]]] - [m_2,[m_5,[m_4,[m_3,m_6]]]]$ \\
$[m_3,[[m_4,m_2],[m_5,m_6]]] = \; $ &  $[m_3,[m_4,[m_2,[m_5,m_6]]]] - [m_3,[m_2,[m_4,[m_5,m_6]]]]$ \\
$[m_3,[[m_2,m_5],[m_4,m_6]]] = \; $ &  $[m_3,[m_2,[m_5,[m_4,m_6]]]] - [m_3,[m_5,[m_2,[m_4,m_6]]]]$ \\
$[m_3,[[m_5,m_4], [m_2,m_6]]] =\; $ &  $[m_3,[m_5,[m_4,[m_2,m_6]]]] - [m_3,[m_4,[m_5,[m_2,m_6]]]]$ \\
$[m_4,[[m_2,m_3],[m_5,m_6]]] = \; $ &  $[m_4,[m_2,[m_3,[m_5,m_6]]]] - [m_4,[m_3,[m_2,[m_5,m_6]]]]$ \\
$[m_4,[[m_5,m_2],[m_3,m_6]]] = \; $ &  $[m_4,[m_5,[m_2,[m_3,m_6]]]] - [m_4,[m_2,[m_5,[m_3,m_6]]]]$ \\
$[m_4,[[m_3,m_5],[m_2,m_6]]] = \; $ &  $[m_4,[m_3,[m_5,[m_2,m_6]]]] - [m_4,[m_5,[m_3,[m_2,m_6]]]]$ \\
\end{tabular}

\renewcommand{\arraystretch}{1.3}
\begin{tabular}{ l l}
$[m_5,[[m_3,m_2],[m_4,m_6]]] = \; $ &  $[m_5,[m_3,[m_2,[m_4,m_6]]]] - [m_5,[m_2,[m_3,[m_4,m_6]]]]$ \\
$[m_5,[[m_2,m_4],[m_3,m_6]]] = \; $ &  $[m_5,[m_2,[m_4,[m_3,m_6]]]] - [m_5,[m_4,[m_2,[m_3,m_6]]]]$ \\
$[m_5,[[m_3,m_4], [m_2,m_6]]]= \; $ &  $[m_5,[m_4,[m_3,[m_2,m_6]]]] - [m_5,[m_3,[m_4,[m_2,m_6]]]]$ \\
\end{tabular}

\renewcommand{\arraystretch}{1.2}
\begin{tabular}{ l r}
$ [m_6,[[m_2,m_3],[m_4,m_5]]]= $ &  $	[m_4,[m_5,[m_3,[m_2,m_6]]]]-[m_5,[m_4,[m_3,[m_2,m_6]]]]$ \\
 & $	-[m_4,[m_5,[m_2,[m_3,m_6]]]]+[m_5,[m_4,[m_2,[m_3,m_6]]]]$\\
& $	+[m_2,[m_3,[m_5,[m_4,m_6]]]]-[m_3,[m_2,[m_5,[m_4,m_6]]]]$\\
& $ -[m_2,[m_3,[m_4,[m_5,m_6]]]]+[m_3,[m_2,[m_4,[m_5,m_6]]]]$ \\
\end{tabular}

\begin{tabular}{ l r}
$[m_6,[[m_4,m_2],[m_3,m_5]]]= $ 
& $	[m_3,[m_5,[m_2,[m_4,m_6]]]]-[m_5,[m_3,[m_2,[m_4,m_6]]]]$ \\
& $	-[m_3,[m_5,[m_4,[m_2,m_6]]]]+[m_5,[m_3,[m_4,[m_2,m_6]]]]$ \\
& $	+[m_4,[m_2,[m_5,[m_3,m_6]]]]-[m_2,[m_4,[m_5,[m_3,m_6]]]]$ \\
& $	-[m_4,[m_2,[m_3,[m_5,m_6]]]]+[m_2,[m_4,[m_3,[m_5,m_6]]]]$ \\
\end{tabular}

\begin{tabular}{ l r}
$[m_6,[[m_2,m_5],[m_3,m_4]]]= $ 
& $	[m_3,[m_4,[m_5,[m_2,m_6]]]]-[m_4,[m_3,[m_5,[m_2,m_6]]]]$ \\
& $	-[m_3,[m_4,[m_2,[m_5,m_6]]]]+[m_4,[m_3,[m_2,[m_5,m_6]]]]$ \\
& $	+[m_2,[m_5,[m_4,[m_3,m_6]]]]-[m_5,[m_2,[m_4,[m_3,m_6]]]]$ \\
& $	-[m_2,[m_5,[m_3,[m_4,m_6]]]]+[m_5,[m_2,[m_3,[m_4,m_6]]]]$ \\
\end{tabular}

\subsection{Proof of proposition \ref{not in image}} \label{proof subsection}
\begin{figure}[ht]
\vspace{.2cm}
\includegraphics[width=6cm]{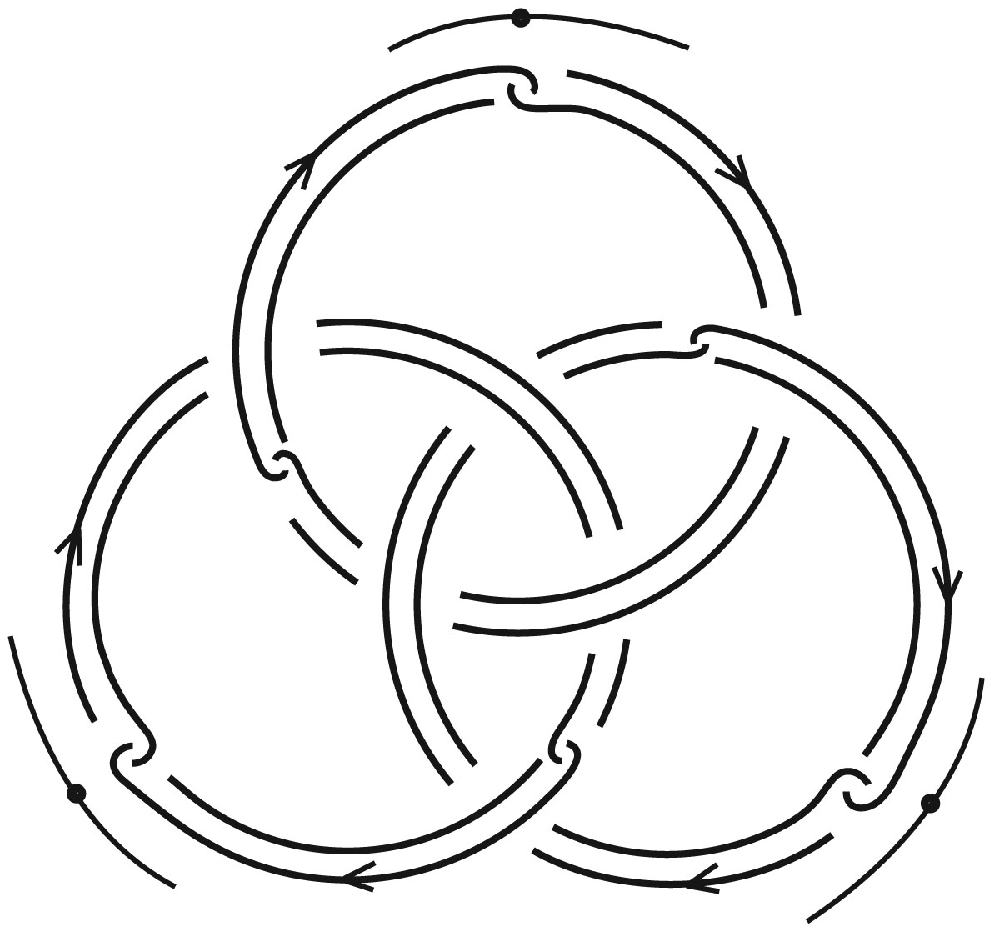}
{\small
    \put(-170,73){$l_6$}
    \put(-135,133){$l_1$}
    \put(-41,133){$l_2$}
    \put(-115,-7){$l_5$}
    \put(-63,-7){$l_4$}
    \put(-3,65){$l_3$}
    \put(-178,30){$r_3$}
    \put(0,30){$r_2$}
    \put(-114,156){$r_1$}}
{\small 
    \put(-300,-7){
\tikzstyle{mybox} = [draw=blue,rectangle, rounded corners, inner sep=2.5pt, inner ysep=2.5pt]
\begin{tikzpicture}
\node [mybox] (box){${\gamma}_1 m_1+{\gamma}_2 m_2+{\gamma}_3 m_3+{\gamma}_4 m_4$};
\end{tikzpicture}}
    \put(-54,148){
\tikzstyle{mybox} = [draw=blue,rectangle, rounded corners, inner sep=2.5pt, inner ysep=2.5pt]
\begin{tikzpicture}
\node [mybox] (box){${\alpha}_3 m_3+{\alpha}_4 m_4+{\alpha}_5 m_5 +{\alpha}_6 m_6$};
\end{tikzpicture}}
    \put(-26,-7){
\tikzstyle{mybox} = [draw=blue,rectangle, rounded corners, inner sep=2.5pt, inner ysep=2.5pt]
\begin{tikzpicture}
\node [mybox] (box){${\beta}_1 m_1+{\beta}_2 m_2+{\beta}_5 m_5+{\beta}_6 m_6$};
\end{tikzpicture}}
}
\vspace{.3cm}
\caption{A hypothetical general solution to the relative-slice
problem in figure \ref{fig:BorDouble}.}
\label{fig:BorDouble2}
\end{figure}

To prove that the link in figure \ref{fig:BorDouble} is not relatively slice up to link-homotopy,
consider a hypothetical solution shown in figure \ref{fig:BorDouble2}.

The circled expressions next to the dotted components show how the slices homologically go over the
attached $2$-handles (only homological information is relevant, since the calculations below
involve $5$-fold commutators and all $6$-fold commutators are trivial in the Milnor group).
Figure \ref{fig:BorDouble2} shows the general solution: according to condition \ref{homological condition} the slices for $l_1,l_2$ do not go over the $2$-handle attached to
$r_1$, and the analogous restriction for the other slices. Here the coefficients ${\alpha}_i, {\beta}_j, {\gamma}_k$ are integers, 
algebraic multiplicities of the hypothetical slices going over the $2$-handles. Interpret this as a band-sum 
of the components $l_1,\ldots, l_6$ with parallel copies of $r_1,r_2,r_3$, yielding a link $(l'_1, \ldots, l'_6)$. 
A direct analogue of proposition \ref{234 homotopically trivial} shows that omitting the first component gives a homotopically trivial link, so $M{\pi}_1(S^3\smallsetminus (l'_2\cup\ldots \cup l'_6))\cong MFree_{m_2,\ldots, m_6}$.
\begin{figure}[h]
\includegraphics[width=8.5cm]{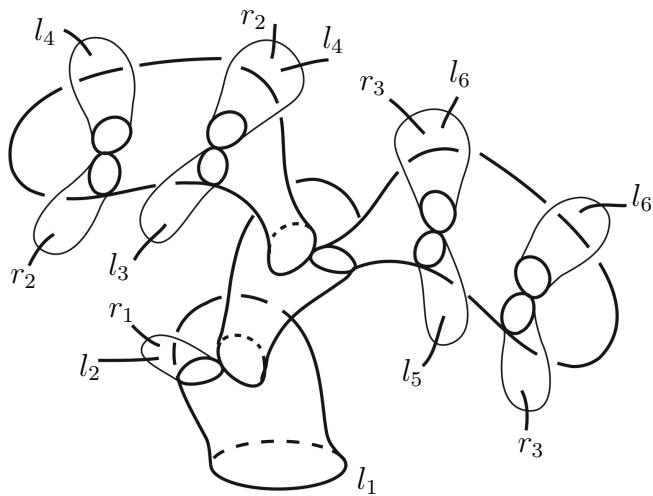}
    \put(-106,3){$l_1$}
    \put(-211,46){$l_2$}
    \put(-200,68){$r_1$}
    \put(-71,157){$l_6$}
    \put(-105,154){$r_3$}
    \put(-2,110){$l_6$}
    \put(-45,17){$r_3$}
    \put(-89,44){$l_5$}
    \put(-150,179){$r_2$}
    \put(-228,173){$l_4$}
    \put(-119,170){$l_4$}
    \put(-200,84){$l_3$}
    \put(-237,81){$r_2$}
\caption{A capped grope bounded by $l_1$ in figure \ref{fig:BorDouble2}.}
\label{fig:Borgrope}
\end{figure}

The component $l_1$ bounds a capped grope in $S^3$, shown in figure \ref{fig:Borgrope}. The body (consisting of the surface stages) of the grope is embedded in the complement of 
the rest of the link in the $3$-sphere, and the caps intersect the components $l_2, \ldots, l_6, r_1, r_2, r_3$ as indicated in the figure. This grope 
is a version of the grope shown in the context of the Hopf link in figures \ref{fig:grope}, \ref{fig:gropelocation}, found in the setting of the Borromean rings. 
It is not a {\em half-grope} considered in section \ref{Hopf subsection}, but rather a grope of a more general type. Specifically, the two third stage surfaces are attached
to a full symplectic basis of the second stage surface, similarly to {\em symmetric gropes} considered in \cite{FQ}. (Since the second stage surface 
has genus one, a symplectic basis consists of two embedded curves forming a basis of the first homology of the surface.) Gropes of the type shown in figure \ref{fig:Borgrope}
are a geometric analogue of commutators of the form (\ref{CommutatorForm}).

%The curve $l'_1$ is a band sum of $l_1$ with ${\beta}_1$ parallel copies of $r_2$ and ${\gamma}_1$ parallel copies of $r_3$. (The curves $r_1, r_2, r_3$ are oriented consistently with the orientations of the components $l_i$, as in figure \ref{fig:gropelocation}. Positive integers ${\beta}_1, {\gamma}_1$ correspond to band sums respecting orientations of the curves $r_2, r_3$ and negative integers correspond to opposite orientations.)  The element representing $l'_1$ in the Milnor group $MFree_{m_2,\ldots, m_6}$ is the sum of the element representing  $l_1$ with ${\beta}_1$ times the element representing $r_2$ and ${\gamma}_1$ times the element representing $r_3$. Here we write ``sum'' since all of them are elements of the fifth term of the lower central series of $MFree_{m_2,\ldots, m_6}$, an abelian group. As shown in figure \ref{fig:Bnew} and discussed just below that figure, each curve $r_2, r_3$ in turn is a band sum of two ``almost parallel'' copies of the relevant components $l_i$. These ``almost parallel'' copies also may be seen as two of the dotted curves in figure \ref{fig:more examples} on the right. Each of them bounds a capped grope almost identical to that in figure \ref{fig:Borgrope}. (The only difference is that the cap for the first stage surface has a single intersection, with the corresponding $l$-component.)

A geometrically transparent way of identifying the individual commutator summands of $l'_1$ of the form  (\ref{CommutatorForm}) is to use the grope splitting technique of
\cite{K6}. Alternatively, this may be done using the commutator identity (\ref{product identity}), generalizing the argument in section \ref{Hopf subsection}.
For the relative slice problem in figure \ref{fig:BorDouble2} to have a link-homotopy solution, $l'_1$ must equal the element $w$ 
in  lemma \ref{intersection prop}. Collecting the coefficients of each commutator in the definition of $w$,
one gets the following system of equations.

%$l_1=[m_2+{\alpha}_3 m_3+{\alpha}_4 m_4+{\alpha}_5 m_5 +{\alpha}_6 m_6, $\\
%$[[m_3+{\beta}_1 m_1+{\beta}_2 m_2+{\beta}_5 m_5+{\beta}_6 m_6, \; \,
%m_4+{\beta}_1 m_1+{\beta}_2 m_2+{\beta}_5 m_5+{\beta}_6 m_6],$\\
%$[m_5+{\gamma}_1 m_1+{\gamma}_2 m_2+{\gamma}_3 m_3+{\gamma}_4 m_4, \; \,
%m_6+{\gamma}_1 m_1+{\gamma}_2 m_2+{\gamma}_3 m_3+{\gamma}_4 m_4]]]$.

\renewcommand{\arraystretch}{1.4}
\begin{tabular}{ l l c }
$[m_2,[[m_3,m_4],[m_5,m_6]]]$ \; \; \; & $1=1$  \; \; \;   & (1) \\
$[m_2,[[m_4,m_5],[m_3,m_6]]]$ \; \; \; & ${\beta}_6{\gamma}_4- {\beta}_5{\gamma}_3=1$  \; \; \; & (2) \\
$[m_2,[[m_5,m_3],[m_4,m_6]]]$ \; \; \;  &  ${\beta}_6{\gamma}_3 - {\beta}_5 {\gamma}_4= 1$  \; \;  \; &  (3) \\
$[m_3,[[m_4,m_2],[m_5,m_6]]]$ \; \; \; & $-{\alpha}_3{\beta}_2+{\beta}_1{\alpha}_4=1$ \; \; \; &  (4) \\
$[m_3,[[m_2,m_5],[m_4,m_6]]]$ \; \; \; & $-{\alpha}_3 {\beta}_6{\gamma}_2-{\beta}_1{\alpha}_5{\gamma}_4=1$ \; \; \; &  (5) \\
$[m_3,[[m_5,m_4],[m_2,m_6]]]$ \; \; \; & ${\alpha}_3{\beta}_5{\gamma}_2 +{\beta}_1{\alpha}_6{\gamma}_4=1$ \; \; \; &  (6) \\
$[m_4,[[m_2,m_3],[m_5,m_6]]]$ \; \; \; & $-{\alpha}_4{\beta}_2+{\beta}_1{\alpha}_3=1$  \; \; \; &  (7) \\
$[m_4,[[m_5,m_2],[m_3,m_6]]]$ \; \; \; &  $-{\alpha}_4 {\beta}_6{\gamma}_2-{\beta}_1{\alpha}_5{\gamma}_3=1$\; \; \; &  (8) \\
$[m_4,[[m_3,m_5],[m_2,m_6]]]$ \; \; \; &  ${\alpha}_4{\beta}_5{\gamma}_2+{\beta}_1{\alpha}_6{\gamma}_3=1$ \; \; \; &  (9) \\
$[m_5,[[m_3,m_2],[m_4,m_6]]]$ \; \; \; & ${\alpha}_5 {\beta}_2{\gamma}_4 -{\gamma}_1{\alpha}_3{\beta}_6  =1$ \; \; \; &  (10) \\
$[m_5,[[m_2,m_4],[m_3,m_6]]]$ \; \; \; &  ${\alpha}_5 {\beta}_2{\gamma}_3-{\gamma}_1{\alpha}_4{\beta}_6    =1$ \; \; \; &  (11) \\
$[m_5,[[m_3,m_4],[m_2,m_6]]]$ \; \; \; & ${\alpha}_5  {\gamma}_2+{\gamma}_1  {\alpha}_6   =1$ \; \; \; &  (12) \\
$[m_6,[[m_2,m_3],[m_4,m_5]]]$ \; \; \; & ${\alpha}_6  {\beta}_2{\gamma}_4-{\gamma}_1 {\alpha}_3{\beta}_5     =1$\; \; \; &  (13) \\
$[m_6,[[m_4,m_2],[m_3,m_5]]]$ \; \; \; & ${\alpha}_6 {\beta}_2{\gamma}_3 -{\gamma}_1{\alpha}_4{\beta}_5    =1$\; \; \; &  (14) \\
$[m_6,[[m_2,m_5],[m_3,m_4]]]$ \; \; \; & ${\alpha}_6 {\gamma}_2 +{\gamma}_1 {\alpha}_5   =1$\; \; \; &  (15) \\
\end{tabular}

The first commutator is already present in the expression for $l_1$ without the relator curves, see (\ref{l6}),
and there are no further contributions from the relator curves, so (1) is automatically satisfied. There are 
14 remaining equations in 12 variables. It is possible to exploit various symmetries of the equations and analyze them by hand,
however for the sake of a concise exposition a Mathematica calculation is included below. There are in fact three families of 
solutions, however there are no {\em integral} solutions that are relevant for the geometric problem at hand. The ``closest'' it gets to 
integers are rational solutions in the first family given in the Mathematica output:

$ {\alpha}_3=-1,  {\alpha}_4=-1,  {\alpha}_5=-2,  {\alpha}_6=-2,$\\
${\beta}_1=1, {\beta}_2=2, {\beta}_5=1, {\beta}_6=-3,$\\
${\gamma}_1=0, {\gamma}_2=-1/2, {\gamma}_3=-1/4, {\gamma}_4=-1/4.$

A non-integral solution as above does not seem to have a geometric interpretation in the context of the relative-slice problem.

\subsection{Mathematica calculation} \label{Mathematica calculation} This section contains the Mathematica
program for solving the system of equations in section \ref{proof subsection}, followed by the program output discussed 
above.

%\begin{doublespace}
\noindent\({\text{Solve}[} \; \; 
%\pmb{b[6]c[4]-b[5]c[3]\text{==}1\&\&}\\
{b[6]c[4]-b[5]c[3]\text{==}1\&\&}\\
{b[6] c[3] - b[5] c[4]\text{==}1\&\&}\\
{-a[3] b[2] + b[1] a[4]\text{==}1\&\&}\\
{ -a[3] b[6] c[2] - b[1] a[5] c[4]\text{==}1\&\&}\\
{a[3] b[5] c[2] + b[1] a[6] c[4]\text{==}1\&\&}\\
{ -a[4] b[2] + b[1] a[3]\text{==}1\&\&}\\
{ -a[4] b[6] c[2] - b[1] a[5] c[3]\text{==}1\&\&}\\
{a[4] b[5] c[2] + b[1] a[6] c[3] == 1\text{==}1\&\&}\\
{a[5] b[2] c[4] - c[1] a[3] b[6]\text{==}1\&\&}\\
{a[5] b[2] c[3] - c[1] a[4] b[6]\text{==}1\&\&}\\
{a[5] c[2] + c[1] a[6]\text{==}1\&\&}\\
{a[6] b[2] c[4] - c[1] a[3] b[5]\text{==}1\&\&}\\
{a[6] b[2] c[3] - c[1] a[4] b[5]\text{==}1\&\&}\\
{a[6] c[2] + c[1] a[5]\text{==}1,}\\
{\{a[3],a[4],a[5],a[6],b[1],b[2],b[5],b[6],c[1],c[2],c[3],c[4]\}]}\)
%\end{doublespace}

%\noindent\(\text{Solve}\text{::}\text{svars}: \text{Equations may not give solutions for all solve variables. }\underline{\rangle\rangle }\)

{\bf Output:}

\noindent\(\left\{\left\{a[3]\to -\frac{1}{b[1]},a[4]\to -\frac{1}{b[1]},a[5]\to -\frac{2 b[5]}{b[1]},a[6]\to -\frac{2 b[5]}{b[1]},b[2]\to 2 b[1],\right.\right.\\
\left.b[6]\to-3 b[5],c[1]\to 0,c[2]\to -\frac{b[1]}{2 b[5]},c[3]\to -\frac{1}{4 b[5]},c[4]\to -\frac{1}{4 b[5]}\right\},\\
\left\{a[3]\to -\frac{\sqrt{3}}{2 b[1]},a[4]\to-\frac{\sqrt{3}}{2 b[1]},a[5]\to 0,a[6]\to \frac{3 \left(-5 b[5]-3 \sqrt{3} b[5]\right)}{2 \left(3 b[1]+2 \sqrt{3} b[1]\right)},\right.\\
b[2]\to \frac{1}{3}\left(3 b[1]+2 \sqrt{3} b[1]\right),b[6]\to -b[5]-\sqrt{3} b[5],c[1]\to -\frac{2 \left(3 b[1]+2 \sqrt{3} b[1]\right)}{3 \left(5+3 \sqrt{3}\right)b[5]},\\
\left.c[2]\to -\frac{2 \left(3 b[1]+2 \sqrt{3} b[1]\right)}{3 \left(5+3 \sqrt{3}\right) b[5]},c[3]\to \frac{-1-\sqrt{3}}{\left(5+3 \sqrt{3}\right)b[5]},c[4]\to \frac{-1-\sqrt{3}}{\left(5+3 \sqrt{3}\right) b[5]}\right\},\\
\left\{a[3]\to \frac{\sqrt{3}}{2 b[1]},a[4]\to \frac{\sqrt{3}}{2 b[1]},a[5]\to0,a[6]\to \frac{3 \left(5 b[5]-3 \sqrt{3} b[5]\right)}{2 \left(-3 b[1]+2 \sqrt{3} b[1]\right)},\right.\\
b[2]\to \frac{1}{3} \left(3 b[1]-2 \sqrt{3} b[1]\right),b[6]\to-b[5]+\sqrt{3} b[5],c[1]\to -\frac{2 \left(-3 b[1]+2 \sqrt{3} b[1]\right)}{3 \left(-5+3 \sqrt{3}\right) b[5]},\\
\left.\left.c[2]\to -\frac{2 \left(-3+2 \sqrt{3}\right)b[1]}{3 \left(-5+3 \sqrt{3}\right) b[5]},c[3]\to \frac{1-\sqrt{3}}{\left(-5+3 \sqrt{3}\right) b[5]},c[4]\to \frac{1-\sqrt{3}}{\left(-5+3 \sqrt{3}\right)
b[5]}\right\}\right\}\)

\bigskip

{\em Acknowledgements}. I would like to thank Michael Freedman for many discussions on the subject, and Jim Conant for sharing his insight on commutator calculus,
in particular an elegant proof of lemma \ref{intersection prop}. I also would like to thank the referee for helpful comments which improved the exposition of the paper.

I am grateful to the Max Planck Institute for Mathematics in Bonn for hospitality and support.


\medskip


\begin{thebibliography}{10}

\bibitem{F2} M. Freedman, {\em A geometric reformulation
of four dimensional surgery}, Topology Appl., 24 (1986), 133-141.

\bibitem{FL} M. Freedman and X.S. Lin, {\em On the $(A,B)$-slice
problem}, Topology Vol. 28 (1989), 91-110.


\bibitem{FK} M. Freedman and V. Krushkal, {\em Topological
arbiters}, J. Topol. 5 (2012), 226-247.

\bibitem{FQ} M. Freedman and F. Quinn, {\em The topology of
4-manifolds}, Princeton Math. Series 39, Princeton, NJ, 1990.

\bibitem{K6} V. Krushkal, {\em Exponential separation in $4$-manifolds}, Geom. Topol. 4 (2000), 397-405.

\bibitem{K4} V. Krushkal, {\em A counterexample to the strong version of Freedman's conjecture},
Ann. of Math. 168 (2008), 675-693.

\bibitem{K5} V. Krushkal, {\em Robust four-manifolds and robust embeddings}, 
Pacific J. Math. 248 (2010), 191-202.

%\bibitem{KT} V. Krushkal and P. Teichner, {\em Alexander duality, gropes and link homotopy}, Geom. Topol. 1 (1997), 51-69.

\bibitem{MKS} W. Magnus, A. Karrass and D. Solitar, Combinatorial group theory: Presentations of groups in terms of generators and relations, Interscience Publishers, New York-London-Sydney 1966.

\bibitem{M} J. Milnor, {\em Link Groups}, Ann. Math 59 (1954), 177-195.

\bibitem{R} C. Reutenauer,  Free Lie algebras, London Mathematical Society Monographs, Oxford University Press, New York, 1993.
\end{thebibliography}
\end{document}